\def\yes{\if00}
\def\no{\if01}
\def\iftwelvept{\yes}
\def\ifusepdf{\no}
\def\ifusepsfont{\no}

\iftwelvept
\documentclass[leqno,12pt]{amsart}
\else
\documentclass[leqno]{amsart}
\fi
\usepackage{amsfonts}
\usepackage{amssymb}
\usepackage{amscd}
\usepackage{latexsym}
\usepackage{verbatim}

\ifusepdf
\usepackage{hyperref}
\else\fi
\ifusepsfont
\usepackage[T1]{fontenc}
\usepackage{times}
\else\fi

\iftwelvept
\else\fi

\newtheorem{Theorem}{Theorem}[section]
\newtheorem{Proposition}[Theorem]{Proposition}
\newtheorem{Lemma}[Theorem]{Lemma}
\newtheorem{Corollary}[Theorem]{Corollary}
\newtheorem{Claim}{Claim}[Theorem]

\theoremstyle{definition}
\newtheorem{Definition}[Theorem]{Definition}
\newtheorem{Remark}[Theorem]{Remark}
\newtheorem{Example}[Theorem]{Example}

\newtheorem{Question}[Theorem]{Question}

\renewcommand{\theTheorem}{\arabic{section}.\arabic{Theorem}}

\newcommand{\ZZ}{{\mathbb{Z}}}
\newcommand{\QQ}{{\mathbb{Q}}}
\newcommand{\RR}{{\mathbb{R}}}
\newcommand{\CC}{{\mathbb{C}}}
\newcommand{\PP}{{\mathbb{P}}}
\newcommand{\OO}{{\mathcal{O}}}

\newcommand{\Spec}{\operatorname{Spec}}
\newcommand{\Pic}{\operatorname{Pic}}
\newcommand{\ord}{\operatorname{ord}}

\newcommand{\Supp}{\operatorname{Supp}}
\newcommand{\zero}{\operatorname{div}}

\newcommand{\calF}{{\mathcal{F}}}
\newcommand{\calX}{{\mathcal{X}}}
\newcommand{\calY}{{\mathcal{Y}}}
\newcommand{\calZ}{{\mathcal{Z}}}
\newcommand{\calW}{{\mathcal{W}}}
\newcommand{\calL}{{\mathcal{L}}}
\newcommand{\calM}{{\mathcal{M}}}
\newcommand{\calB}{{\mathcal{B}}}
\newcommand{\calH}{{\mathcal{H}}}
\newcommand{\calJ}{{\mathcal{J}}}
\newcommand{\calN}{{\mathcal{N}}}

\newcommand{\bd}[1]{{\pmb{#1}}}

\newcommand{\Rea}{\operatorname{Re}}
\newcommand{\rank}{\operatorname{rk}}
\newcommand{\vol}{\operatorname{vol}}
\newcommand{\cherncl}{\operatorname{c}}
\newcommand{\acherncl}{\operatorname{\widehat{c}}}
\newcommand{\adeg}{\operatorname{\widehat{deg}}}
\newcommand{\Gal}{\operatorname{Gal}}
\newcommand{\wlim}{\operatorname{w-lim}}

\newcommand{\canmet}{\widehat{\Vert\cdot\Vert}_{L, \bd{\overline{f}}}}
\newcommand{\canmetS}{\widehat{\Vert\cdot\Vert}_{L, S(\bd{\overline{f}})}}
\newcommand{\canmetSi}{\widehat{\Vert\cdot\Vert}_{L, S^i(\bd{\overline{f}})}}
\newcommand{\canmetP}{\widehat{\Vert\cdot\Vert}_{\OO_{\PP^N}(1), \bd{\overline{f}}}}

\newcommand{\Proof}{{\sl Proof.}\quad}

\newcommand{\QED}{{\unskip\nobreak\hfil\penalty50\quad\null\nobreak\hfil
{$\Box$}\parfillskip0pt\finalhyphendemerits0\par\medskip}}
\newcommand{\rest}[2]{\left.{#1}\right\vert_{{#2}}}
\begin{document}

%\begin{flushright}
%{\large \fbox{1st Draft}}
%\end{flushright}
%\bigskip

%%%%%%%%%%%
%% Title %%
%%%%%%%%%%%
\title[Canonical heights for random iterations]%
{Canonical heights for random iterations in certain varieties}
\author{Shu Kawaguchi}
\address{Department of Mathematics, Faculty of Science,
Kyoto University, Kyoto, 606-8502, Japan}
\email{kawaguch@math.kyoto-u.ac.jp}
%\date{\DateTime}
%\date{\DateTime, (\Version)}
%\thanks{{\em 2000 Mathematics Subject Classification:} 11G50}
\begin{abstract}
We show the existence of canonical heights (normalized heights) 
of subvarieties for bounded sequences of morphisms 
and give some applications. 
\end{abstract} 

\maketitle
\section*{Introduction}
\renewcommand{\theTheorem}{\Alph{Theorem}} 
Theory of heights of points (or subvarieties in general) 
is a basic tool in Diophantine geometry.  
Particular heights that enjoy nice properties, called ``canonical
heights'' or ``normalized heights,'' are sometimes of great use as are
N\'eron-Tate heights on abelian varieties (see, for example, 
Hindry--Silverman \cite[Part~B, Part~F.2]{HS}).
On the other hand, 
dynamics of sequences of polynomial mappings (random iterations) 
in $\PP^N$ over $\CC$ have been studied in Russakovskii--Sodin 
\cite{RS} and Forn{\ae}ss--Weickert \cite{FW} among others.  
In particular, in \cite{FW} invariant currents that transform 
well relative to such sequences and their average have been considered. 

In this paper, we show the existence of heights 
of subvarieties that transform well relative to 
``bounded'' sequences of such morphisms and give some applications. 
Although our interest mainly lies in morphisms of $\PP^N$ 
of degree greater than or equal to $2$, our results 
also hold in the setting of Call--Silverman \cite{CS} and Zhang 
\cite{Zh}: Namely, let $X$ be a projective variety over a number field $K$ 
and $L$ a line bundle on $X$, and we consider morphisms $f: X \to X$ 
such that $f^*{L} \simeq L^{\otimes d_f}$ for some integer 
$d_f \geq 2$. 

Let $\bd{f}=(f_i)_{i=1}^{\infty}$ be a sequence of morphisms such 
that $f_i: X \to X$ satisfies $f_i^*{L} \simeq L^{\otimes d_{f_i}}$ 
for some integer $d_{f_i} \geq 2$ for each $i$, and we consider 
the iterations by $\bd{f}$: 
\[
X \overset{f_1}{\longrightarrow} X 
\overset{f_2}{\longrightarrow} X 
\overset{f_3}{\longrightarrow} \cdots
\]

Let $h_L$ be a height function on $X(\overline{K})$ 
corresponding to $L$. 
Then we know $c(f_i) := \sup_{x \in X(\overline{K})}
\left\vert \frac{1}{d_{f_i}}h_L(f_i(x)) - h_L(x) 
\right\vert < + \infty$ for each $i$. 
We say that a sequence $\bd{f}$ is {\em bounded} if 
$c(\bd{f}) := \sup_{i\geq 1} c(f_i) < + \infty$. 
Notice that, although we use $h_L$ to define boundedness, 
the property that $\bd{f}$ is bounded is independent of  
the choice of height functions $h_L$ corresponding to $L$. 
Let $S$ be the shift map which sends 
$\bd{f}=(f_i)_{i=1}^{\infty}$ to 
$S(\bd{f}) := (f_{i+1})_{i=1}^{\infty}$. 
See \S\ref{sec:canonical:height:functions} for more details. 

After reviewing some basic facts on heights in \S\ref{sec:preliminaries}, 
we show the following theorem in \S\ref{sec:canonical:height:functions}. 

\begin{Theorem}[cf. Theorem~\ref{thm:can:height:function}]
\label{thm:A}
There is a unique way to attach to each bounded sequence 
$\bd{f}=(f_i)_{i=1}^{\infty}$  
a height function 
$\widehat{h}_{L, \bd{f}}: X(\overline{K}) \to \RR$ 
such that 
\begin{enumerate}
\item[(i)] $\sup_{x \in X(\overline{K})} 
\left\vert \widehat{h}_{L, \bd{f}}(x) - h_L(x)\right\vert
\leq 2 c(\bd{f})$\textup{;} 
\item[(ii)] $\widehat{h}_{L, S(\bd{f})} \circ f_1 
= d_{f_1}\; \widehat{h}_{L, \bd{f}}$. 
\end{enumerate}
The function $\widehat{h}_{L, \bd{f}}$ is independent of the choice 
of height functions $h_L$ corresponding $L$.

Moreover, if $L$ is ample, then $\widehat{h}_{L, \bd{f}}$ is 
non-negative, and $\widehat{h}_{L, \bd{f}}(x) = 0$ if and only if 
$x$ is $\bd{f}$-preperiodic. 
\textup{(}Here a point $x \in X(\overline{K})$ is said to 
be $\bd{f}$-preperiodic if the forward 
orbit of $x$ under $\bd{f}$, i.e., 
$O^+_{\bd{f}}(x) := \{ x,\; f_1(x),\; f_2(f_1(x)),\; f_3(f_2 (f_1(x))),\; 
\cdots \}$ is finite.\textup{)} 

We call $\widehat{h}_{L, \bd{f}}$ a canonical height function 
\textup{(}normalized height function\textup{)} for $\bd{f}$. 
\end{Theorem}

Note that, for a single morphism $f$, 
i.e., when $f_1 = f_2 = \cdots (=f)$, 
a canonical height function for $f$ is constructed in \cite{CS} and \cite{Zh}. 
We remark that some conditions (such as ``boundedness'') on $\bd{f}$ are
necessary to ensure the existence of height functions that behave 
well relative to $\bd{f}$. Indeed, 
dropping the assumption of boundedness, there 
is {\em not} in general a way to attach to each 
sequence $\bd{f}=(f_i)_{i=1}^{\infty}$ a height function 
$\widehat{h}_{L, \bd{f}}$ such that (i) $\widehat{h}_{L, \bd{f}} = h_L + O(1)$ 
and (ii) $\widehat{h}_{L, S(\bd{f})} \circ f_1 
= d_{f_1}\; \widehat{h}_{L, \bd{f}}$ (cf. Example~\ref{eg:3}). 

Iterations by finitely many morphisms as follows 
give examples of bounded sequences: 
Let $g_1, \ldots, g_k: X \to X$ be
morphisms such that $g_j^*{L} \simeq L^{\otimes d_{g_j}}$ 
for some integer $d_{g_j} \geq 2$ for $j = 1, \ldots, k$. 
We set $J := \{1, \ldots, k\}$ and $W := \prod_{i=1}^{\infty} J$. 
For $w = (w_i)_{i=1}^{\infty} \in W$, we set 
$\bd{f}_{w} = (g_{w_i})_{i=1}^{\infty}$. 
Then $\bd{f}_{w}$ is a bounded sequence.  
As an immediate consequence of Theorem~\ref{thm:A}, 
we get the following corollary due to  
Masseron \cite[\S2.9]{Ma}. 

\begin{Corollary}[\cite{Ma}, 
see also Corollary~\ref{cor:can:height:function}(1)(2)]
\label{cor:B}
Assume $L$ is ample. 
Then for any positive integer $D$, the set 
\[
\left\{
x \in X(\overline{K}) 
\;\left\vert\;
\begin{gathered}
\text{$[K(x):K] \leq D$}, \\
\text{$x$ is $\bd{f}_w$-preperiodic 
for some $w \in W$}  
\end{gathered}
\right. \right\}
\]
is finite.  
\end{Corollary}

In \S\ref{sec:averaging}, 
noting that averaging of currents over 
a space with a suitable probability measure and topology is studied 
in \cite{FW}, we consider averaging of canonical height functions 
$\widehat{h}_{L, \bd{f}_w}$ where $w$ runs over $W$. 
By \cite[Theorem~1.2]{Ka}, 
there is a unique height function 
$\widehat{h}_{L, \{g_1, \ldots, g_k\}}$ on $X(\overline{K})$
with the following two properties: 
(i) $\widehat{h}_{L, \{g_1, \ldots, g_k\}} = h_L + O(1)$; and 
(ii) $\sum_{i=1}^k \widehat{h}_{L, \{g_1, \ldots, g_k\}} \circ g_i 
= (d_{g_1} + \cdots + d_{g_k}) \;\widehat{h}_{L, \{g_1, \ldots, g_k\}}$.
Then we have the following proposition. 

\begin{Proposition}[cf. Proposition~\ref{prop:averaging}]
\label{prop:C}
We give $J$ the discrete topology, 
and let $\nu$ be the measure on $J$ that assigns 
mass $\frac{d_{g_j}}{d_{g_1} + \cdots + d_{g_k}}$ to $j \in J$. 
Let $\mu := \prod_{i=1}^{\infty} \nu$ be the product measure 
on $W$. Then we have, for all $x \in X(\overline{K})$, 
$\widehat{h}_{L, \{g_1, \ldots, g_k\}}(x) 
= \int_W \widehat{h}_{L, \bd{f}_{w}}(x) d\mu(w)$. 
\end{Proposition}

In \S\ref{sec:canonical:heighs:for:subvarieties}, 
we consider canonical heights of subvarieties of $X_{\overline{K}}$ 
for $\bd{f} = (f_i)_{i=1}^{\infty}$, as in \cite{Zh} for 
a single morphism $f$.  
We  show in Theorem~\ref{thm:can:height} that, 
assuming some bounded conditions on $\bd{f}$, 
there exist a height $\widehat{h}_{L, \bd{f}}(Y)$
of any subvariety $Y$ of $X_{\overline{K}}$ that 
behaves well relative to $\bd{f}$ as in Theorem~\ref{thm:A}. 
(However, the statement on preperiodicity in Theorem~\ref{thm:A} 
is changed to the following: If 
$\{Y,\; f_1(Y),\; f_2(f_1(Y)),\; \cdots \}$ is a finite set, then 
$\widehat{h}_{L, \bd{f}}(Y) = 0$.) 
In \S\ref{sec:canonical:heights:for:subvarieties:2}, 
we give another construction of canonical heights of subvarieties 
by using adelic sequences.  

In \S\ref{sec:admissible:metrics}, we consider a local theory 
in the setting of \cite[\S2]{Zh}. 
Let $X$ a projective variety over an algebraically 
closed valuation field $K_v$, 
and $L$ a line bundle on $X$. Let $\Vert\cdot\Vert$ be 
a continuous and bounded metric on $L$. 
We consider a pair $\overline{f} := (f, \varphi)$ 
such that $f: X \to X$ is a morphism over $K_v$ and 
$\varphi: L^{\otimes d_{f}} \overset{\sim}{\rightarrow} f^*L$ 
for some integer $d_f \geq 2$. Then 
we will show in Theorem~\ref{thm:admissible:metric} that, 
assuming some boundedness 
conditions on $\bd{\overline{f}} = (\overline{f_i})_{i=1}^{\infty}$, 
there exists a bounded and continuous metric 
$\canmet$ on $L$ satisfying $\varphi_1^* f_1^* \canmetS
= \canmet^{d_{f_1}}$, where $S$ is the shift map. 
We call $\canmet$ the admissible metric for 
$\bd{\overline{f}}$. 

Suppose that $K_v=\CC$, $X=\PP^N$ and 
$L=\OO_{\PP^N}(1)$ is equipped with the Fubini-Study metric 
$\Vert\cdot\Vert_{FS}$ in the above setting.  
Let $\calJ_{\bd{\overline{f}}}$ be the smallest closed 
set on $\PP^N(\CC)$ such that on its complement the family 
$\{f_1,\; f_2\circ f_1,\; f_3\circ f_2\circ f_1,\;\cdots\}$ is normal. 
As constructed in \cite{FW}, one has 
the Green current $T_{\bd{\overline{f}}}$ on $\PP^N(\CC)$, 
whose support coincides with $\calJ_{\bd{\overline{f}}}$ 
(Proposition~\ref{prop:estimate:Green:function}). 
Then, as noted in \cite[Proposition~3.3.1]{Ka} for a single morphism 
$f$, we find $\cherncl_1(\OO_{\PP^N}(1), \canmetP) 
= T_{\bd{\overline{f}}}$, where the left-hand side is the first Chern 
current of $\OO_{\PP^N}(1)$ with the admissible metric 
$\canmetP$ (Proposition~\ref{prop:Chern:current}). 

Finally, on $\PP^1$, 
we show equidistribution of small points 
for random iterations as follows. Let $K$ be a number field, 
and we fix an embedding $\overline{K} \hookrightarrow \CC$ 
so that we regard $\PP^1(\overline{K})$ as a subset of $\PP^1(\CC)$. 
Let $\bd{\overline{f}} = \left((f_i, \varphi_i)\right)_{i=1}^{\infty}$ 
be a sequence of $K$-morphisms on $\PP^1$ satisfying some 
boundedness conditions (see \S\ref{sec:equidistribution}). 
Let $\widehat{h}_{\OO_{\PP^1}(1), \bd{f}}: 
\PP^1(\overline{K}) \to \RR$ be the 
canonical height function for $\bd{f} = (f_i)_{i=1}^{\infty}$, 
and $T_{\bd{\overline{f}}}$ the Green current for $\bd{\overline{f}}$ 
on $\PP^1(\CC)$. 
For a point $x \in \PP^1(\overline{K})$, we denote the Galois orbit 
of $x$ by $G(x) := \{g(x) \mid g \in \Gal(\overline{K}/K)\}$. 
Then we have the following theorem.

\begin{Theorem}[cf. Theorem~\ref{thm:equidistribution}]
\label{thm:D}
Let $(x_j)_{j=1}^{\infty}$ be a sequence of distinct points 
in $\PP^1(\overline{K})$ such that 
$\lim_{j\to +\infty} \widehat{h}_{\OO_{\PP^1}(1), 
\bd{f}}(x_j) = 0$. Then 
$\wlim_{j\to +\infty} \frac{1}{\# G(x_j)} \sum_{y \in G(x_j)} \delta_{y}
= T_{\bd{\overline{f}}}$. 
\textup{(}Here $\delta_{y}$ denotes the Dirac measure of mass $1$ 
supported in $y$.\textup{)}
\end{Theorem}

In \cite{SUZ}, Szpiro--Ullmo--Zhang proved equidistribution of small
points on abelian varieties, and ever since much progress has been
seen around this subject. On $\PP^1$, equidistribution of small points
will follow from Autissier \cite[Proposition~4.1.4 and Remarque in
Introduction]{Au}.  Baker--Hsia \cite{BH} explicitly stated and proved
equidistribution of small points on $\PP^1$ for polynomial maps $f$
over global fields satisfying the product formula.  See also
Chamber-Loir \cite{CL}, Baker--Rumely \cite{BR} and
Favre--Rivera-Letelier \cite{FR} for equidistribution of small points
on $\PP^1$ over archimedean and non-archimedean fields.  
Our proof of Theorem~\ref{thm:D} uses Arakelov geometry. 
Somewhat different ingredients are that we
use an estimate of analytic torsions of Bismut--Vasserot \cite{BV} and
Vojta \cite{Vo} and that we use only intersections of $C^{\infty}$ hermitian
line bundles (i.e., not using $L_1^2$-intersection theory), 
in hoping that this might be useful to consider a higher
dimensional case (cf. Question~\ref{quesiton:analytic:torsion}).

\smallskip
{\sl Acknowledgment.}\quad 
My deep thanks are due to Professor Laurent Denis. Indeed, 
this paper is motivated by e-mail from him in October, 2004. He kindly gave me valuable comments and suggestions, and informed me of Dr.\ Masseron's work \cite{Ma}. I would also like to thank Professors Atsushi Moriwaki and Joseph Silverman for valuable comments. 

\renewcommand{\theTheorem}{\arabic{section}.\arabic{Theorem}}
\setcounter{equation}{0}
\section{Preliminaries}
\label{sec:preliminaries}
In this section, we briefly review some basic facts on heights of
closed points and subvarieties. 

Let $h_{nv}: \PP^N({\QQ})\to\RR$ be the logarithmic naive height
function: For a number field $K$ and $x=(x_0:\cdots:x_N) \in \PP^{N}(K)$,  
$h_{nv}(x)$ is given by
\[
h_{nv}(x) = 
\frac{1}{[K:\QQ]}
\left[
\sum_{P\in\Spec(O_K)\setminus\{0\}} 
\max_{i}\{\log\Vert x_i\Vert_P\} + 
\sum_{\sigma: K \hookrightarrow \CC} 
\max_{i}\{\log|\sigma(x_i)|\}
\right],
\]
where $O_K$ is the ring of integers of $K$ and $\Vert x_i\Vert_P =
\#(O_K/P)^{-\ord_P(x_i)}$.

We will need the following two theorems on height functions. 
For details of height functions, we refer to \cite{HS}. 

\begin{Theorem}[Weil's height machine, \cite{HS} Theorem~B.3.2]
\label{thm:properties:of:height:functions}
There is a way to attach to any projective variety $X$ 
over $\overline{\QQ}$
and any line bundle $L$ on $X$ a function
\[
h_{X,L}: X(\overline{\QQ}) \to \RR
\]
with the following three properties\textup{:}
\begin{enumerate}
\item[(i)]
$h_{X,L\otimes M} = h_{X,L} + h_{X,M} + O(1)$ 
for any line bundles $L$ and $M$ on $X$\textup{;} 
\textup{(}Here $O(1)$ denotes a bounded function 
on $X(\overline{\QQ})$\textup{;)}
\item[(ii)]
If $X = \PP^N$ and $L = \OO_{\PP^N}(1)$, then 
$h_{\PP^N, \OO_{\PP^N}(1)} = h_{nv} + O(1)$\textup{;}
\item[(iii)]
If $f: Y \to X$ is a morphism of projective varieties 
and $L$ is a line bundle on $X$, then $h_{Y,f^*L} = h_{X,L}\circ f + O(1)$. 
\end{enumerate}
Moreover, height functions $h_{X,L}$ are determined up to $O(1)$
by the above three properties.
\end{Theorem}

We say that $h_{X,L}$ is a height function (or a Weil height function) 
corresponding to $L$. We often write $h_L$ for $h_{X,L}$ 
when $X$ is clear from the context.

\begin{Theorem}[\cite{HS} Theorem~B.2.3, Theorem~B.3.2]
\label{thm:height:function:2}
Assume $L$ is ample. Let $h_{X,L}$ be a height function corresponding to $L$. 
\begin{enumerate}
\item[(1)] \textup{(}Northcott's finiteness theorem\textup{)} 
For any real number $c$ and positive integer $D$, the set 
\[
\{x \in X(\overline{\QQ}) \mid 
[\QQ(x):\QQ] \leq D, h_{X,L}(x) \leq c \}
\]
is finite. 
\item[(2)] \textup{(}positivity\textup{)} 
There is a constant $c'$ such that 
$h_{X,L}(x) \geq c'$ for all $x \in X(\overline{\QQ})$. 
\end{enumerate}
\end{Theorem} 

Next we recall some properties of heights of subvarieties. 
We refer to \cite{SABK} for details of arithmetic intersection theory, 
and to \cite{BGS} for details of heights of subvarieties. 

Let $K$ be a number field, and $O_K$ its ring of integers. 
Let $X$ be a projective variety over $K$, and $L$ a line bundle on $X$.   
We say that $(\calX, \overline{\calL})$ is a {\em $C^{\infty}$ model} of $(X,L)$ if 
$\calX$ is a projective arithmetic variety over $O_K$ (i.e., 
an integral scheme projective and flat over $\Spec(O_K)$) 
that extends $X$ and if 
$\overline{\calL}$ is a $C^{\infty}$ hermitian $\QQ$-line bundle on $\calX$ that 
extends $L$. Here $\overline{\calL}$ is called a $C^{\infty}$ hermitian $\QQ$-line 
bundle if $\overline{\calL}$ is 
a pair $(\calL, \{\Vert\cdot\Vert_{\sigma}\}_{\sigma})$ 
such that $\calL$ is a $\QQ$-line bundle on $\calX$ 
and $\Vert\cdot\Vert_{\sigma}$'s are $C^{\infty}$ hermitian metrics on 
$\calL \otimes_{K^{\sigma}}\CC$ for the embeddings
$\sigma: K \hookrightarrow \CC$ 
which are invariant under complex conjugation. 

In the rest of this section, we assume $L$ is ample. 
We fix a $C^{\infty}$ model $(\calX, \overline{\calL})$ of $(X,L)$. 

Let $Y$ be a subvariety of $X_{\overline{K}}$. Take a finite 
extension field $K'$ of $K$ such that $Y$ is defined over 
$K'$. Let $O_{K'}$ be the ring of integers of $K'$, 
and $p: \calX\times_{\Spec(O_K)}\Spec(O_K') \to \calX$ 
the natural morphism. Let $\calY$ be the Zariski closure 
of $Y$ in $\calX\times_{\Spec(O_K)}\Spec(O_K')$. 
Then the height of $Y$ with respect to $(\calX, \overline{\calL})$ 
is defined by 
\begin{equation}
\label{eqn:height}
h_{(\calX, \overline{\calL})}(Y) 
: = \frac{\adeg(\acherncl_1(\rest{p^*\overline{\calL}}{\calY})^{\dim Y +1})}%
{[K':\QQ] (\dim Y +1) \deg( \rest{L}{Y}^{\cdot\dim Y})} \quad \in\RR. 
\end{equation}

\begin{Theorem}%
[\cite{BGS} Proposition~3.2.2]
\label{thm:two:models}
Let $(\calX_1, \overline{\calL_1})$ and 
$(\calX_2, \overline{\calL_2})$ be two $C^{\infty}$ models of $(X,L)$, where  
we assume $L$ is ample.
Then there exists a constant $C$ such that for any subvariety 
$Y$ of $X_{\overline{K}}$, one has 
\[
\left\vert
h_{(\calX_1, \overline{\calL_1})}(Y) - 
h_{(\calX_2, \overline{\calL_2})}(Y)
\right\vert \leq C.
\]
\end{Theorem}

\setcounter{equation}{0}
\section{Canonical height functions}
\label{sec:canonical:height:functions}
As in \S\ref{sec:preliminaries}, let $K$ be a number field, 
$X$ a projective variety over $K$, and $L$ 
a line bundle on $X$. 
We fix a height function $h_{L}: X(\overline{K})\to\RR$ corresponding 
to $L$. 

Let $\calH$ be the set of all morphisms $f: X \to X$ over 
$K$ such that $f^*(L) \simeq L^{\otimes d_f}$ for some integer 
$d_f \geq 2$. 
For $f\in \calH$, we set
\[
c(f) := \sup_{x \in X(\overline{K})} 
\left\vert \frac{1}{d_{f}} h_L(f(x)) - h_L(x) \right\vert. 
\]
Since $\frac{1}{d_{f}} h_L\circ f - h_L$ is a bounded function on
$X(\overline{K})$ by Theorem~\ref{thm:properties:of:height:functions}(iii), 
we see that $c(f) < +\infty$.

Let $\bd{f}=(f_i)_{i=1}^{\infty}$ 
be a sequence with $f_i \in \calH$
for $i \geq 1$. The set of all such sequences is denoted by
$\prod_{i=1}^{\infty} \calH$. For $\bd{f} \in \prod_{i=1}^{\infty} \calH$, 
we set
\[
c(\bd{f}) := \sup_{i \geq 1} c(f_i) \; \in \RR \cup \{+\infty\}. 
\]
When $c(\bd{f}) <+\infty$, 
we say that $\bd{f}$ is a {\em bounded} sequence. 

Suppose $h_L^{\prime}$ be another height function corresponding
to $L$. Then there is a constant $c$ such that 
$\sup_{x \in X(\overline{K})} 
\left\vert h_L(x) - h_L^{\prime}(x) \right\vert \leq c$, 
and thus 
\[
\sup_{x \in X(\overline{K})} 
\left\vert \frac{1}{d_{f}} h_L^{\prime}(f(x)) - h_L^{\prime}(x) 
\right\vert 
\leq \sup_{x \in X(\overline{K})} 
\left\vert \frac{1}{d_{f}} h_L(f(x)) - h_L(x) 
\right\vert  + \left(\frac{1}{d_f} + 1\right) c 
\leq c(f) + 2 c.
\] 
This shows that the property that $\bd{f}$ is bounded 
is independent of the choice of height functions 
$h_L$ corresponding to $L$. 

Let $\calB$ be the set of 
all bounded sequences in $\prod_{i=1}^{\infty} \calH$. 
For any nonnegative number $c$, 
we define the subset $\calB_{c}$ of $\calB$ by 
\[
\calB_{c} := 
\{\bd{f} =(f_i)_{i=1}^{\infty} \in \calB 
\mid c(\bd{f}) \leq c \}.
\]

Let $S: \prod_{i=1}^{\infty} \calH \to \prod_{i=1}^{\infty} \calH$ be 
the {\em shift map} which sends 
$\bd{f} =(f_i)_{i=1}^{\infty} \in\prod_{i=1}^{\infty} \calH$ to 
$S(\bd{f}) := (f_{i+1})_{i=1}^{\infty} \in\prod_{i=1}^{\infty} \calH$. 
Since $c(S(\bd{f})) \leq c(\bd{f})$, $S$ maps $\calB$ into $\calB$, 
and $\calB_{c}$ into $\calB_{c}$ for any $c$. 

\medskip
For $\bd{f} =(f_i)_{i=1}^{\infty} \in \prod_{i=1}^{\infty} \calH$ 
and $x \in X(\overline{K})$,  
we consider iterated points by $\bd{f}$ with the initial point $x$:  
\[
x \mapsto f_1(x) \mapsto f_2(f_1(x)) \mapsto f_3(f_2(f_1(x))) 
\mapsto \cdots .
\]

The set $\{ x, f_1(x),\; f_2(f_1(x)),\; f_3(f_2 (f_1(x))),\; 
\cdots \}$ is
called the {\em forward orbit} of $x$ under $\bd{f}$, and is denoted
by $O^+_{\bd{f}}(x)$. A point $x \in X(\overline{K})$ is said to be
{\em $\bd{f}$-preperiodic} if $O^+_{\bd{f}}(x)$ is a finite set.
Note that, when $f_1 = f_2 = \cdots (= f)$, the forward orbit 
under $\bd{f}$ is just the forward orbit under $f$ in the usual sense, 
and an $\bd{f}$-preperiodic point is just a usual $f$-preperiodic 
point. 

\begin{Example}
\label{eg:1}
We consider iterations by a finite number of morphisms 
$g_1, \ldots, g_k \in \calH$. 
We set $J:= \{1, \ldots, k\}$ and 
$W := \prod_{i=1}^{\infty} J$. 
For $w = (w_i)_{i=1}^{\infty} \in W$, 
we set $\bd{f}_{w} = (g_{w_i})_{i=1}^{\infty}$. 
Put $c := \max\{c(g_1), \ldots, c(g_k)\}$. Then 
$\left\{\bd{f}_{w} \;\left\vert\; w \in W \right. \right\}
\subset \calB_c$.  
\end{Example}

\begin{Example}
\label{eg:2}
We give some concrete examples. Let $K=\QQ$, $X=\PP^N$ and 
$L=\OO_{\PP^N}(1)$. For a fixed height function, we take the 
naive height function $h_{nv}: \PP^N(\overline{\QQ}) \to\RR$. 
\begin{enumerate}
\item[(1)]
For each $m \geq 2$, let 
$g_m^{\prime}: \PP^N\to\PP^N$ be the morphism of degree 
$d_{g_m^{\prime}}=m$ defined by 
\[
g_m^{\prime}(x_0:\cdots:x_N) = (x_0^m:\cdots:x_N^m).
\]
Then $c(g_m^{\prime}) = 0$. Thus any sequence 
$\bd{f} = (f_i)_{i=1}^{\infty}$ such that, for each $i$, 
there is $m_i$ with $f_i = g_{m_i}^{\prime}$ belongs to $\calB_0$.   
\item[(2)] 
For each $m \geq 2$, let 
$g_m^{\prime\prime}: \PP^N\to\PP^N$ be the morphism 
of degree $d_{g_m^{\prime\prime}}=m$ defined by 
\[
g_m^{\prime\prime}(x_0:x_1:\cdots:x_N) = (x_0^m+x_1^m:x_1^m:\cdots:x_N^m).
\] 
Then it is easy to see $c(g_m^{\prime\prime}) \leq \log 2$. Thus any sequence 
$\bd{f} = (f_i)_{i=1}^{\infty}$ such that, for each $i$, 
there is $m_i$ with $f_i = g_{m_i}^{\prime\prime}$ 
belongs to $\calB_{\log 2}$.   
\end{enumerate}
\end{Example}

\begin{Theorem}
\label{thm:can:height:function}
Let $X$ be a projective variety over a number field $K$, and $L$ 
a line bundle on $X$. Let 
$h_{L}: X(\overline{K})\to\RR$ be a height function 
corresponding to $L$. 
\begin{enumerate}
\item[(1)]
There is a unique way to attach to each bounded sequence 
$\bd{f}=(f_i)_{i=1}^{\infty} \in \calB$  
a height function 
\[
\widehat{h}_{L, \bd{f}}: X(\overline{K}) \to \RR
\]
such that 
\begin{enumerate}
\item[(i)] $\sup_{x \in X(\overline{K})} 
\left\vert \widehat{h}_{L, \bd{f}}(x) - h_L(x)\right\vert
\leq 2 c(\bd{f})$\textup{;} 
\item[(ii)] $\widehat{h}_{L, S(\bd{f})} \circ f_1 
= d_{f_1}\; \widehat{h}_{L, \bd{f}}$. 
\end{enumerate}
Moreover, $\widehat{h}_{L, \bd{f}}$ is independent of the choice of 
height functions $h_L$ corresponding to $L$. 
\item[(2)]
Assume $L$ is ample. Then $\widehat{h}_{L, \bd{f}}$ satisfies 
the following properties\textup{:} 
\begin{enumerate}
\item[(iii)] $\widehat{h}_{L, \bd{f}}(x) \geq 0$ for 
all $x \in X(\overline{K})$\textup{;} 
\item[(iv)] $\widehat{h}_{L, \bd{f}}(x) = 0$ if and only if 
$x$ is $\bd{f}$-preperiodic.
\end{enumerate}
\end{enumerate}
We call $\widehat{h}_{L, \bd{f}}$ a canonical height function 
\textup{(}normalized height function\textup{)} for $\bd{f}$. 
\end{Theorem}

\Proof 
(1) We first construct $\widehat{h}_{L, \bd{f}}$ for $\bd{f} \in
\calB$. By the definition of $c(\bd{f})$, we have for any $i$ 
\[
\sup_{x \in X(\overline{K})} \left\vert
\frac{1}{d_{f_i}} h_{L}(f_i(x)) - h_L(x) 
\right\vert \leq c(\bd{f}) \quad (< +\infty).
\] 
For $i=0$ we set $h_0 := h_L$, and for $i \geq 1$ we set 
\[
h_i := \frac{1}{\prod_{\alpha=1}^i d_{f_{\alpha}}} 
h_L \circ f_i \circ f_{i-1} \circ\cdots\circ f_1.
\]

\begin{Claim}
For $x \in X(\overline{K})$, 
$\{h_i(x)\}_{i=0}^{\infty}$ is a Cauchy sequence. 
\end{Claim}

Indeed, we have 
\begin{align*}
\left\vert {h}_{i+1}(x) - h_i(x) \right\vert
& = 
\left\vert \frac{1}{\prod_{\alpha=1}^{i+1} d_{f_{\alpha}}}{h}_{L}
(f_{i+1}\circ \cdots \circ f_1(x)) - 
\frac{1}{\prod_{\alpha=1}^{i}d_{f_{\alpha}}}{h}_{L}
(f_{i}\circ \cdots \circ f_1(x)) \right\vert \\
& = 
\frac{1}{\prod_{\alpha=1}^{i} d_{f_{\alpha}}}
\left\vert \frac{1}{d_{f_{i+1}}} 
{h}_{L}(f_{i+1} (f_i \circ \cdots \circ f_1(x))) - 
{h}_{L}(f_{i}\circ \cdots \circ f_1(x)) \right\vert \\
& \leq 
\frac{c(\bd{f})}{\prod_{\alpha=1}^{i} d_{f_{\alpha}}}
\leq 
\frac{c(\bd{f})}{2^{i}}.
\end{align*}
Thus we get the claim. We define $\widehat{h}_{L, \bd{f}}$ to be 
\begin{equation}
\label{eqn:can:height:fcn}
\widehat{h}_{L, \bd{f}}(x)
:= \lim_{i\to\infty} h_i(x)
= \lim_{i\to\infty} \frac{1}{\prod_{\alpha=1}^{i}d_{f_{\alpha}}}
{h}_{L}(f_{i}\circ \cdots \circ f_1(x)).
\end{equation}

Let us check that $\widehat{h}_{L, \bd{f}}(x)$ satisfies the conditions
(i) and (ii). Since
\begin{equation}
\label{eqn:can:height:fcn:2}
\left\vert {h}_{i}(x) - h_L(x) \right\vert
\leq 
\sum_{\alpha=0}^{i-1} 
\left\vert {h}_{\alpha+1}(x) - h_{\alpha}(x) \right\vert
\leq 
\sum_{\alpha=0}^{i-1}\frac{c(\bd{f})}{2^{\alpha}}
\leq 2 c(\bd{f}),
\end{equation}
we obtain (i) by letting $i$ to the infinity. 

Since $\widehat{h}_{L, S(\bd{f})}(y) = 
\lim_{i \to\infty}
\frac{1}{\prod_{\alpha=1}^i
  d_{f_{\alpha+1}}} h_L(f_{i+1} \circ f_{i} \circ\cdots\circ f_2(y))$,
substituting $y=f_1(x)$ gives 
\begin{align*}
\widehat{h}_{L, S(\bd{f})}(f_1(x)) 
& = 
\lim_{i\to\infty}
\frac{1}{\prod_{\alpha=1}^i d_{f_{\alpha+1}}} 
h_L(f_{i+1} \circ f_{i} \circ\cdots\circ f_2 \circ f_1(x)) \\
& = 
d_{f_1} 
\lim_{i\to\infty}
\frac{1}{\prod_{\alpha=1}^{i+1} d_{f_{\alpha}}} 
h_L(f_{i+1} \circ f_{i} \circ\cdots\circ f_2 \circ f_1(x)) 
= 
d_{f_1} \widehat{h}_{L, \bd{f}}(x). 
\end{align*}
Thus we get (ii). 

To show the uniqueness of 
$\{\widehat{h}_{L, \bd{f}}\}_{\bd{f}\in\calB}$, suppose
$\{\widehat{h}_{L, \bd{f}}^{\prime}\}_{\bd{f}\in\calB}$ 
are functions satisfying
(i) and (ii). By (ii), we have
\begin{align*}
\sup_{x \in X(\overline{K})} \left\vert
\widehat{h}_{L, \bd{f}}^{\prime}(x) - \widehat{h}_{L, \bd{f}}(x)
\right\vert 
& \leq 
\frac{1}{d_{f_1}} 
\sup_{x \in X(\overline{K})} \left\vert
\widehat{h}_{L, S(\bd{f})}^{\prime}(x) - \widehat{h}_{L, S(\bd{f})}(x)
\right\vert \\
& \leq \cdots \leq
\frac{1}{\prod_{\alpha=1}^{i}d_{f_{\alpha}}} 
\sup_{x \in X(\overline{K})} \left\vert
\widehat{h}_{L, S^i(\bd{f})}^{\prime}(x) - \widehat{h}_{L, S^i(\bd{f})}(x)
\right\vert \leq 
\frac{4 c(S^{i}(\bd{f}))}{\prod_{\alpha=1}^{i}d_{f_{\alpha}}}. 
\end{align*}

Since $c(S^{i}(\bd{f})) \leq c(\bd{f})$, we obtain 
$\sup_{x \in X(\overline{K})} \left\vert
\widehat{h}_{L, \bd{f}}^{\prime}(x) - \widehat{h}_{L, \bd{f}}(x)
\right\vert = 0$ by letting $i$ to the infinity, whence 
$\widehat{h}_{L, \bd{f}} = \widehat{h}_{L, \bd{f}}^{\prime}$. 

Since the difference of two height functions corresponding to $L$ 
is bounded on $X(\overline{K})$, 
it follows from \eqref{eqn:can:height:fcn} that 
$\widehat{h}_{L, \bd{f}}$ is independent of the choice of 
height functions $h_L$. 

(2) Now assuming $L$ is ample, let us see (iii) and (iv).  Since $h_L$
is bounded below by Theorem~\ref{thm:height:function:2}(2), we get
(iii) by \eqref{eqn:can:height:fcn}.

Next we show (iv). Suppose $\widehat{h}_{L, \bd{f}}(x) =0$. Take a finite
extension field $K'$ of $K$ such that $x$ is defined over $K'$. Then
for any $i \geq 1$, $f_i\circ\cdots\circ f_1(x)$ is also defined over
$K'$. We set
\[
T := 
\{y \in X(K') \mid h_L(y) \leq 2 c(\bd{f})\}.
\]
We claim that the forward orbit $O^+_{\bd{f}}(x)$ is contained in $T$.
Indeed, since $\widehat{h}_{L, \bd{f}}(x) =0$ and $\left\vert
  \widehat{h}_{L, \bd{f}}(x) - h_L(x) \right\vert \leq 2 c(\bd{f})$, we see
$x\in T$. For $i\geq 1$, we have
\[
\widehat{h}_{L, S^i(\bd{f})}(
f_i\circ\cdots\circ f_1(x)) = 
\left( \prod_{\alpha=1}^i d_{f_{\alpha}} \right)
\widehat{h}_{L, \bd{f}}(x) =0.
\]
It follows from  
\begin{equation}
\label{eqn:canonical:height:function:3}
\left\vert \widehat{h}_{L, S^i(\bd{f})}(f_i\circ\cdots\circ f_1(x)) -
  h_L(f_i\circ\cdots\circ f_1(x)) \right\vert \leq 2 c(S^i(\bd{f})) \leq
2 c({\bd{f}})
\end{equation}
that $h_L(f_i\circ\cdots\circ f_1(x)) \leq 2 c(\bd{f})$. Thus we get the
claim. Since $T$ is a finite set by Northcott's finiteness theorem
(Theorem~\ref{thm:height:function:2}), so is $O^+_{\bd{f}}(x)$. In
other words, $x$ is $\bd{f}$-preperiodic.

Finally we will show that $x$ is $\bd{f}$-preperiodic then
$\widehat{h}_{L, \bd{f}}(x) =0$. To see this, suppose
$\widehat{h}_{L, \bd{f}}(x) =a > 0$. Then
\[
\widehat{h}_{L, S^i(\bd{f})}(
f_i\circ\cdots\circ f_1(x)) = 
\left(\prod_{\alpha=1}^i d_{f_{\alpha}}\right)
a \to +\infty
\]
as $i$ tends to the infinity. 
Thus by \eqref{eqn:canonical:height:function:3}, 
$O^+_{\bd{f}}(x)$ cannot be finite. 
\QED

The following corollary is an immediate consequence 
of Theorem~\ref{thm:can:height:function}. 
We remark that Corollary~\ref{cor:can:height:function}(2) 
is due to Masseron \cite[\S2.9]{Ma}. 

\begin{Corollary}
\label{cor:can:height:function}
Assume $L$ is ample. 
\begin{enumerate}
\item[(1)]
Let $c$ be a nonnegative number, and $D$ a positive integer. Then 
the set 
\[
\bigcup_{\bd{f}\in\calB_c} 
\left\{
x \in X(\overline{K}) 
\mid 
\text{$[K(x):K] \leq D$, $x$ is $\bd{f}$-preperiodic}  
\right\}
\]
is finite. 
\item[(2)]
Let $g_1, \ldots, g_k$ be elements of $\calH$. 
As in Example~\ref{eg:1}, we set 
$W := \prod_{i=1}^{\infty} \{1, \ldots, k\}$, and 
$\bd{f}_{w} = (g_{w_i})_{i=1}^{\infty}$ for 
$w=(w_i)_{i=1}^{\infty} \in W$. 
Then for any positive integer $D$, the set 
\[
\left\{
x \in X(\overline{K}) 
\;\left\vert\;
\begin{gathered}
\text{$[K(x):K] \leq D$}, \\
\text{$x$ is $\bd{f}_w$-preperiodic 
for some $w \in W$}  
\end{gathered}
\right. \right\}
\]
is finite.  
\end{enumerate}
\end{Corollary}

\Proof
By Theorem~\ref{thm:can:height:function}, we obtain for $\bd{f}\in\calB_c$, 
%\begin{multline*}
\[
\left\{x \in X(\overline{K}) 
\;\left\vert\;
\begin{gathered}
\text{$[K(x):K] \leq D$}, \\
\text{$x$ is $\bd{f}$-preperiodic} 
\end{gathered}
\right.\right\}  
= 
\left\{
x \in X(\overline{K}) 
\;\left\vert\;
\begin{gathered}
\text{$[K(x):K] \leq D$}, \\ 
\text{$\widehat{h}_{L,\bd{f}}(x) = 0$}   
\end{gathered}
\right.\right\}. 
\]
%\end{multline*}
Since 
$\sup_{x \in X(\overline{K})} \left\vert
\widehat{h}_{L, \bd{f}}(x) - h_L(x) \right\vert \leq 2 c(\bd{f}) \leq 2 c$, 
we find that the above set is a subset of 
\[
\{x \in X(\overline{K}) 
\mid
[K(x):K] \leq D, \; 
{h}_L(x) \leq 2 c \}.
\]

This set is finite by Northcott's finiteness theorem
(Theorem~\ref{thm:height:function:2}), 
and is independent of $\bd{f} \in \calB_c$. Thus
we get the assertion (1). 
Since $\{\bd{f}_{w} \mid w \in W\}$ is contained in 
$\calB_c$ with $c := \max\{c(g_1), \ldots, c(g_k)\}$, 
the assertion (2) follows from (1).  
\QED

\begin{Example}
\label{eg:3}
By dropping the assumption of boundedness in 
Theorem~\ref{thm:can:height:function}, 
one can ask the following question: 
Is there a way to attach to each sequence 
$\bd{f} \in \prod_{i=1}^{\infty}\calH$  
a height function $\widehat{h}_{L, \bd{f}}: X(\overline{K}) \to \RR$ 
such that (i) $\widehat{h}_{L, \bd{f}} = h_L + O(1)$ 
and (ii) $\widehat{h}_{L, S(\bd{f})} \circ f_1 
= d_{f_1}\; \widehat{h}_{L, \bd{f}}$? 

The following example shows that this is not possible in general, 
and thus some conditions (such as ``boundedness'') on $\bd{f}$ are necessary 
to ensure the existence of height functions that behave well relative 
to $\bd{f}$. 

We define $F_i$ inductively: 
\begin{align*}
F_1(x) & = x(x-1), \\ 
F_i(x) & = x \left(
x - F_{i-1}\circ\cdots\circ F_1(i) \right) \quad (i \geq 2).
\end{align*} 
We define $f_i: \PP^1 \to \PP^1$ by 
\[
f_i((x_0:x_1)) = 
\left(x_0^2 : F_i\left(\frac{x_1}{x_0}\right) x_0^2\right). 
\]

Then $f_i$ is defined over $\QQ$ and $d_{f_i} = 2$ for any $i \geq 1$. 
We put $p_i = (1:i) \in \PP^1(\QQ)$ for $i \geq 0$.  
Then $f_{\alpha}(p_0) = p_0$ for any $\alpha \geq 1$, 
and $f_{i}\circ f_{i-1}\circ\cdots\circ f_1(p_i) = p_0$. 

Suppose there exists a height function 
$\widehat{h}_{\OO_{\PP^1}(1), \bd{f}}: \PP^1(\overline{\QQ})\to\RR$ 
satisfying (i) $\widehat{h}_{\OO_{\PP^1}(1), \bd{f}} = h_{nv} + O(1)$ and 
(ii) $\widehat{h}_{\OO_{\PP^1}(1), S^i(\bd{f})} \circ f_1 
= 2 \widehat{h}_{\OO_{\PP^1}(1), S^{i-1}(\bd{f})}$ for $i \geq 1$. 
Then by (ii), 
\begin{align*}
\widehat{h}_{\OO_{\PP^1}(1), \bd{f}} (p_i) 
& = \frac{1}{2^i} 
\widehat{h}_{\OO_{\PP^1}(1), S^i(\bd{f})}( f_i \circ f_{i-1}\circ\cdots\circ f_1(p_i))  \\
& =  \frac{1}{2^i} 
\widehat{h}_{\OO_{\PP^1}(1), S^i(\bd{f})}(p_0) = 
\frac{1}{2^i} 
\widehat{h}_{\OO_{\PP^1}(1), S^i(\bd{f})}( f_i \circ f_{i-1}\circ\cdots\circ f_1(p_0)) 
= \widehat{h}_{\OO_{\PP^1}(1), \bd{f}} (p_0) 
\end{align*}
for any $i \geq 1$. Thus the set 
$T := \{ x \in \PP^1(\QQ) \mid  
\widehat{h}_{\OO_{\PP^1}(1), \bd{f}}(x) 
= \widehat{h}_{\OO_{\PP^1}(1), \bd{f}}(p_0)\}$ 
is an infinite set. This contradicts (i), 
as $\widehat{h}_{\OO_{\PP^1}(1), \bd{f}}$ 
does not satisfy Northcott's finiteness property. 
\end{Example}

We also give two remarks on possible extentions of Theorem~\ref{thm:can:height:function}. 

\begin{Remark}
\label{rmk:R:bundle}
We replace an (integral) line bundle $L$ by an $\RR$-line bundle, and integers $d_{f_i}$ by real numbers. We fix a positive real number $\varepsilon >0$. 
Let $X$ be a projective variety over a number field $K$, and $L$ an 
$\RR$-line bundle on $X$. We consider sequences 
$\bd{f}=(f_i)_{i=1}^{\infty}$, where $f_i: X \to X$ is a morphism such that 
$f_i^*(L) \simeq L^{\otimes d_{f_i}}$ for some real number 
$d_{f_i} \geq (1+\varepsilon)$. Then, with the same definition of boundedness, the corresponding statements to Theorem~\ref{thm:can:height:function} hold. Note that the property of Theorem~\ref{thm:can:height:function} (1)(i) is replaced by 
\begin{enumerate}
\item[(i)']
$\sup_{x \in X(\overline{K})} 
\left\vert \widehat{h}_{L, \bd{f}}(x) - h_L(x)\right\vert
\leq \left(1 + \frac{1}{\varepsilon}\right) c(\bd{f})$.
\end{enumerate}
\end{Remark}

\begin{Remark}
\label{rmk:overline:Q}
We replace a number field $K$ by $\overline{\QQ}$. Let $X$ be a projective variety over $\overline{\QQ}$, and $L$ a
line bundle on $X$. We consider sequences 
$\bd{f}=(f_i)_{i=1}^{\infty}$, where $f_i: X \to X$ is a morphism 
over $\overline{\QQ}$ such that $f_i^*(L) \simeq L^{\otimes d_{f_i}}$ for some integer $d_{f_i} \geq 2$. Then, with the same definition of boundedness, the corresponding statements to Theorem~\ref{thm:can:height:function}(1) and (2)(i) hold. However, the  corresponding statement to Theorem~\ref{thm:can:height:function}(2)(ii) does not hold, and must be replaced by 
\begin{enumerate}
\item[(iv)']
$\widehat{h}_{L, \bd{f}}(x) = 0$ if $x$ is $\bd{f}$-preperiodic. 
\end{enumerate}

To see that the ``only if" part of (iv)' fails in general, we give an example as follows. Let $p_1, p_2, \ldots$ be distinct prime numbers. Let $\zeta_i$ be a $p_i$-th primitive root of unity for each $i$. We set $X = \PP^1_{\overline{\QQ}}$ and $L = \OO_{\PP^1_{\overline{\QQ}}}(1)$, and consider the logarithmic naive height function $h_{nv}$. We define $f_i: \PP^1_{\overline{\QQ}} \to \PP^1_{\overline{\QQ}}$ by 
$f_i((x_0:x_1)) = (\zeta_i x_0^2:x_1^2)$. 
Then
\[
h_{nv}\left(f_i((x_0:x_1))\right) = h_{nv}\left((\zeta_i x_0^2:x_1^2)\right) 
= 2 h_{nv}\left((x_0:x_1)\right). 
\]
This shows that $\bd{f} = (f_i)_{i=1}^{\infty}$ is bounded, and 
$\widehat{h}_{L, \bd{f}} = h_{nv}$. Now we take $x_0 = (1:1)$. 
On one hand, $\widehat{h}_{L, \bd{f}}(x_0) = h_{nv}(x_0) =0$. 
On the other hand, 
since $f_i \circ \cdots \circ f_1 (x_0) = (\zeta_i \zeta_{i-1}^2 \cdots \zeta_{1}^{2^i} : 1)$, $O_{\bd{f}}^+(x_0)$ is an infinite set. Thus 
the ``only if" part of (iv)' fails in this example. 
\end{Remark}

\setcounter{equation}{0}
\section{Averaging of canonical height functions}
\label{sec:averaging}
Let $X$ be a projective variety over a number field $K$, 
and $L$ a line bundle on $X$. 
As in Example~\ref{eg:1}, let $g_1, \ldots, g_k$ be elements of 
$\calH$, and we set $J := \{1, \ldots, k\}$, $W:= \prod_{i=1}^{\infty} J$, 
and $\bd{f}_{w} = (g_{w_i})_{i=1}^{\infty}$ for 
$w = (w_i)_{i=1}^{\infty} \in W$. 
Let $\widehat{h}_{L, \bd{f}_w}$ be the canonical height function 
for $\bd{f}_{w}$. 

In this section, 
noting that averaging of currents over a space with 
a suitable probability measure and topology is studied 
in \cite{FW}, we would like to consider averaging 
of $\widehat{h}_{L, \bd{f}_w}$ over $W$. 

First we recall \cite[Theorem~1.2]{Ka}. 
Since $g_j^*L \simeq L^{\otimes d_{g_j}}$, we have 
$g_1^*L\otimes\cdots\otimes g_k^*L \simeq 
L^{\otimes (d_{g_1} + \cdots + d_{g_k})}$. Thus $(X; g_1, \ldots, g_k)$ 
becomes a particular case of what we call a dynamical eigensystem 
for $L$ of degree $d_{g_1} + \cdots + d_{g_k}$. 
Then we have the canonical height function 
\[
\widehat{h}_{L, \{g_1, \ldots, g_k\}}: X(\overline{K}) \to \RR
\]
for $(X; g_1, \ldots, g_k)$ characterized by the following 
two properties: 
(i) $\widehat{h}_{L, \{g_1, \ldots, g_k\}} = h_L + O(1)$; and 
(ii) $\sum_{j=1}^k \widehat{h}_{L, \{g_1, \ldots, g_k\}} \circ g_j 
= (d_{g_1} + \cdots + d_{g_k}) \;\widehat{h}_{L, \{g_1, \ldots, g_k\}}$.

\begin{Proposition}
\label{prop:averaging}
We give $J$ the discrete topology, 
and let $\nu$ be the measure on $J$ that assigns 
mass $\frac{d_{g_j}}{d_{g_1} + \cdots + d_{g_k}}$ to $j \in J$. 
Let $\mu := \prod_{i=1}^{\infty} \nu$ be the product measure 
on $W$. Then we have, for all $x \in X(\overline{K})$, 
\[
\widehat{h}_{L, \{g_1, \ldots, g_k\}}(x) 
= \int_W \widehat{h}_{\bd{f}_{w}}(x) d\mu(w). 
\]
\end{Proposition}

\Proof
As in Example~\ref{eg:1}, we set 
$c : = \max\{c(g_1), \ldots, c(g_k)\}$. 
For $w = (w_i)_{i=1}^{\infty} \in W$ and $x \in X(\overline{K})$, 
we put 
\[
h_{i}(x, w) : = \frac{1}{\prod_{\alpha=1}^i d_{g_{w_{\alpha}}}}
h_L(g_{w_i}\circ\cdots\circ g_{w_1}(x)).
\]
By \eqref{eqn:can:height:fcn}, 
$\lim_{i \to \infty} h_i(x,w) = \widehat{h}_{L, \bd{f}_{w}}(x)$. 
By \eqref{eqn:can:height:fcn:2}, 
$\vert h_{i}(x, w) \vert \leq 2c +  \vert h_L(x) \vert$. 
Then, if we fix $x$, 
$\{h_i(x, w)\}_{i=1}^{\infty}$ is 
a sequence of integral functions on $W$ bounded by 
$2 c + \vert h_L(x) \vert$. The Lebesgue convergence theorem
implies 
\begin{equation}
\label{eqn:averaging:1}
\int_W \widehat{h}_{L, \bd{f}_{w}}(x) d\mu(w) 
= \lim_{i\to\infty} \int_W h_{i}(x, w) d\mu(w).  
\end{equation}

We set $\widehat{h}^{\prime}(x) 
:= \int_W \widehat{h}_{L, \bd{f}_{w}}(x) d\mu(w)$. 
We will show that as a function on $X(\overline{K})$, 
$\widehat{h}^{\prime}$ satisfies the properties (i) 
and (ii) of $\widehat{h}_{L, \{g_1, \ldots, g_k\}}$

Since 
\[
\left\vert \widehat{h}^{\prime}(x) - h_L(x) \right\vert 
\leq 
\left\vert \int_W\left(
\widehat{h}_{L, \bd{f}_{w}}(x) - h_L(x)
\right) d\mu(w) \right\vert \leq 2 c, 
\]
we get (i). 
To show (ii), we note the following equality: 
\begin{multline*}
\int_W h_{i}(x, w) d\mu(w) = 
\int_W \frac{h_L(g_{w_i}\circ\cdots\circ g_{w_1}(x))}%
{\prod_{\alpha=1}^i d_{g_{w_{\alpha}}}} d\mu(w) \\
= 
\sum_{w_1, \ldots, w_i = 1}^k 
\frac{h_L(g_{w_i}\circ\cdots\circ g_{w_1}(x))}%
{\prod_{\alpha=1}^i d_{g_{w_{\alpha}}}}
\prod_{\alpha=1}^i \frac{d_{g_{w_{\alpha}}}}{d_{g_1}+ \cdots + d_{g_k}}
=  
\sum_{w_1, \ldots, w_i = 1}^k 
\frac{h_L(g_{w_i}\circ\cdots\circ g_{w_1}(x))}%
{(d_{g_1}+ \cdots + d_{g_k})^i}.
\end{multline*}
Then we have 
\begin{align*}
\sum_{j=1}^{k} \widehat{h}^{\prime}(g_j(x)) 
& = \lim_{i \to\infty} \sum_{j=1}^k 
\int_W h_i (g_j(x), w) d\mu(w) \\
& = \lim_{i\to\infty} \sum_{j=1}^{k} \sum_{w_1, \ldots, w_i=1}^k  
\frac{h_L(g_{w_i}\circ\cdots\circ g_{w_1}(g_j(x)))}%
{(d_{g_1}+ \cdots + d_{g_k})^i} \\
& = \lim_{i \to\infty} (d_{g_1}+ \cdots + d_{g_k}) 
\int_W h_{i+1}(x, w) d\mu(w) 
= (d_{g_1}+ \cdots + d_{g_k})  \widehat{h}^{\prime}(x).  
\end{align*}
Thus we get (ii). By the uniqueness of 
$\widehat{h}_{L, \{g_1, \ldots, g_k\}}$, we get the assertion. 
\QED 

\begin{Remark}
As in \cite[Remark~1]{FW}, it is possible to consider other index sets 
$J$ with suitable probability measures and topologies.
\end{Remark}

\setcounter{equation}{0}
\section{Canonical heights of subvarieties, I}
\label{sec:canonical:heighs:for:subvarieties}
As in \S\ref{sec:preliminaries}, let $K$ be a number field, 
$O_K$ its ring of integers, $X$ a projective variety 
over $K$, and $L$ a line bundle on $X$. We assume $L$ is ample. 
We fix a $C^{\infty}$ model $(\calX, \overline{\calL})$ of $(X,L)$. 

Let $f$ be an element of $\calH$, that is,  
$f: X \to X$ is a morphism over $K$ such that 
$f^*(L) \simeq L^{\otimes d_f}$ for some integer $d_f \geq 2$. 
Since we assume $L$ is ample, $f$ is finite. 

\begin{Lemma}
\label{lemma:c:tilde}
There exists a constant $\widetilde{c}$ such that, 
for any subvariety $Y$ of $X_{\overline{K}}$, one has 
\[
\left\vert
\frac{1}{d_f}
h_{(\calX, \overline{\calL})}(f(Y)) 
- h_{(\calX, \overline{\calL})}(Y) 
\right\vert \leq \widetilde{c}. 
\]
\end{Lemma}

\Proof
Let $\calX_1$ be the normalization of the morphism 
$X \overset{f}{\longrightarrow} X \hookrightarrow \calX$, 
and  $\widetilde{f}: \calX_1 \to \calX$ the induced morphism. 
Then $(\calX_1, \widetilde{f}^*(\overline{\calL}))$ is a $C^{\infty}$ model of 
$(X, L^{\otimes d_f})$. 
By Theorem~\ref{thm:two:models}, there exists 
a constant $\tilde{c}^{\prime}$ such that, 
for any subvariety $Y$ of $X_{\overline{K}}$, 
\[
\left\vert
h_{(\calX_1, \widetilde{f}^*(\overline{\calL}))}(Y) 
- h_{(\calX, \overline{\calL}^{\otimes d_f})}(Y) 
\right\vert \leq \tilde{c}^{\prime}. 
\]
On the other hand, by the projection formula, we have  
$h_{(\calX_1, \widetilde{f}^*(\overline{\calL}))}(Y) 
= h_{(\calX, \overline{\calL})}(f(Y))$. 
Since $h_{(\calX, \overline{\calL}^{\otimes d_f})}(Y)  = 
d_f h_{(\calX, \overline{\calL})}(Y)$, 
we obtain the assertion with 
$\widetilde{c}= \frac{\tilde{c}^{\prime}}{d_f}$. 
\QED

For $f \in \calH$, we set 
\[
\widetilde{c}(f) := \sup_{Y \subset X_{\overline{K}}} 
\left\vert
\frac{1}{d_f}
h_{(\calX, \overline{\calL})}(f(Y)) 
- h_{(\calX, \overline{\calL})}(Y) 
\right\vert \;\; \in \RR,  
\]
where $Y$ runs the set of all subvarieties of $X_{\overline{K}}$. 
For $\bd{f} = (f_i)_{i=1}^{\infty} \in \prod_{i=1}^{\infty}\calH$, 
we set 
\[
\widetilde{c}(\bd{f}) := \sup_{i\geq 1} 
\widetilde{c}(f_i) \;\; \in \RR\cup \{+\infty\}.  
\]
We define 
$\widetilde{\calB} := 
\{
\bd{f} = (f_i)_{i=1}^{\infty} \in \prod_{i=1}^{\infty}\calH
\mid \widetilde{c}(\bd{f}) < +\infty \}$.  
Note that, by virtue of Theorem~\ref{thm:two:models}, 
the property that $\bd{f}$ belongs to $\widetilde{\calB}$ 
is independent of the choice of $C^{\infty}$ models $(\calX, \overline{\calL})$ 
of $(X,L)$. 

\begin{Theorem}
\label{thm:can:height}
Let $X$ a projective variety 
over a number field $K$, and $L$ an ample line bundle on $X$. 
Let $(\calX, \overline{\calL})$ be a $C^{\infty}$ model of $(X,L)$. 
\begin{enumerate}
\item[(1)]
There is a unique way to attach to each sequence 
$\bd{f}=(f_i)_{i=1}^{\infty} \in \widetilde{\calB}$ 
\[
\widehat{h}_{L, \bd{f}}: 
\{\text{subvariety of $X_{\overline{K}}$}\} \to \RR 
\]
with the following properties\textup{:}
\begin{enumerate}
\item[(i)]
For any subvariety $Y$ of $X_{\overline{K}}$, one has 
$\left\vert \widehat{h}_{L, \bd{f}}(Y) - 
h_{(\calX, \overline{\calL})}(Y) \right\vert 
\leq 2 \widetilde{c}(\bd{f})$\textup{;} 
\item[(ii)]
$\widehat{h}_{L, S(\bd{f})} \circ f_1 = d_{f_1} \widehat{h}_{L, \bd{f}}$. 
\end{enumerate}
\item[(2)]
Moreover, $\widehat{h}_{L, \bd{f}}$ is independent of the choice of 
$C^{\infty}$ models $(\calX, \overline{\calL})$ of $(X,L)$, and  
satisfies the following properties\textup{:} 
\begin{enumerate}
\item[(iii)]
$\widehat{h}_{L, \bd{f}}(Y) \geq 0$ for any subvariety $Y$ 
of $X_{\overline{K}}$\textup{;}  
\item[(iv)] 
If $\{Y, f_1(Y), f_2(f_1(Y)), \cdots \}$ is a finite set, then 
$\widehat{h}_{L, \bd{f}}(Y) = 0$\textup{;} 
\item[(v)]
If $Y$ is a closed point, then 
$\widehat{h}_{L, \bd{f}}(Y)$ coincides with the one constructed in 
Theorem~\ref{thm:can:height:function}.
\end{enumerate}
\end{enumerate}
We call $\widehat{h}_{L, \bd{f}}(Y)$ a canonical height 
\textup{(}normalized height\textup{)} of $Y$ for $\bd{f}$. 
\end{Theorem}

\Proof
Since we can prove Theorem~\ref{thm:can:height}(1) just  
as in Theorem~\ref{thm:can:height:function}, we only sketch a 
proof here. To construct $\widehat{h}_{L, \bd{f}}(Y)$, we put 
\begin{equation}
\label{eqn:can:height}
h_i(Y) := 
\frac{1}{\prod_{\alpha=1}^i d_{f_{\alpha}}} 
h_{(\calX, \overline{\calL})}(f_i\circ f_{i-1}\circ \cdots \circ f_1(Y)). 
\end{equation}
Then $\{h_i(Y)\}_{i=0}^{\infty}$ is a Cauchy sequence, 
which allows one to define 
$\widehat{h}_{L, \bd{f}}(Y) = \lim_{i\to\infty} h_i(Y)$. 
One can check that $\widehat{h}_{L, \bd{f}}$ satisfies (i) and (ii), 
and is unique. 

We will show (2). 
By virtue of Theorem~\ref{thm:two:models}, it follows 
\eqref{eqn:can:height} that $\widehat{h}_{L, \bd{f}}$ 
is independent of the choice of $C^{\infty}$ models of $(X, L)$. 
We take a $C^{\infty}$ model $(\calX, \overline{\calL})$ of $(X, L)$ 
such that $\calL$ is ample, $\cherncl_1(\overline{L}_{\sigma})$ 
is semipositive on $X_{\sigma}(\CC)$ for each $\sigma: K \to \CC$,  
and $H^0(\calX, \calL^{\otimes n})$ is generated by 
$\{s \in H^0(\calX, \calL^{\otimes n}) \mid \Vert s\Vert_{\sup} < 
1\}$ (cf. \cite[Lemma~1.3]{Mor}). 
Then, for any $i \geq 1$,  $h_i(Y) \geq 0$ for this model. 
Then we get (iii) by letting $i$ to the infinity. 
The finiteness of 
$\{Y, f_1(Y), f_2(f_1(Y)), \cdots \}$ implies 
that of 
$\{ h_{(\calX, \overline{\calL})}(Y),\; 
h_{(\calX, \overline{\calL})}(f_1(Y)),\; 
h_{(\calX, \overline{\calL})}(f_2(f_1(Y))),\; \cdots\}$. 
By (i) and (ii), we have 
\[
\left\vert
d_{f_1} \cdots d_{f_i} \widehat{h}_{L, \bd{f}}(Y)
- h_{(\calX, \overline{\calL})}(f_i \circ \cdots \circ f_1(Y))
\right\vert
\leq
2 \widetilde{c}\left(S^i({\bd{f}})\right)
\leq 
2 \widetilde{c}(\bd{f}) 
\]
for any $i \geq 1$. Thus we get (iv). 
By \eqref{eqn:can:height}, we get (v). 
\QED

\setcounter{equation}{0}
\section{Canonical heights for subvarieties, II}
\label{sec:canonical:heights:for:subvarieties:2}
In this section, we give another construction of 
a canonical height of subvarieties, using the adelic intersection 
theory in Zhang \cite{Zh}. 
This construction will be used 
in \S\ref{sec:equidistribution}. 

Let $\calX$ be a projective arithmetic variety and 
$\overline{\calL} = (\calL, \{\Vert\cdot\Vert_{\sigma}\}_{\sigma})$ 
a $C^{\infty}$ hermitian $\QQ$-line bundle. Following \cite[\S1]{Mor}, 
we say that $\overline{\calL}$ is {\em nef} if 
$\cherncl_1(\calL \otimes_{K^{\sigma}} \CC)$ is semipositive 
on $\calX_{\sigma}(\CC)$ for all embeddings $\sigma: K \hookrightarrow 
\CC$, and $\adeg(\rest{\overline{\calL}}{\Gamma}) \geq 0$ 
for all one-dimensional integral closed subschemes 
$\Gamma$ of $\calX$. 
We say that $\overline{\calL}$ is {\em $\QQ$-effective} if 
there is a positive integer $n$ and a non-zero section 
$s \in H^0(\calX, \calL^{\otimes n})$ such that 
$\Vert s \Vert_{\sigma, \sup} \leq 1$ for all 
$\sigma: K \hookrightarrow \CC$. We write $\overline{\calL} 
\succsim 0$ if $\overline{\calL}$ is $\QQ$-effective. 
Moreover, if $U$ is a non-empty Zariski open set of 
$\calX$ with $\zero(s) \subseteq \calX\setminus U$, then 
we write $\overline{\calL} \succsim_U 0$. 
For two $C^{\infty}$ hermitian $\QQ$-line bundles $\overline{\calL}$ 
and $\overline{\calM}$ on $\calX$, we write 
$\overline{\calL} \succsim \overline{\calM}$  
(resp. $\overline{\calL} \succsim_{U} \overline{\calM}$) 
if $\overline{\calL} \otimes \overline{\calM}^{\otimes -1} 
\succsim 0$ (resp. $\overline{\calL} \otimes \overline{\calM}^{\otimes -1} 
\succsim_U 0$).

In the rest of this section, 
let  $X$ be a projective variety over a number field $K$, 
and $L$ an ample line bundle on $X$.  

Following \cite[\S1]{Zh} and \cite[\S3.1]{Mor}, we define 
an adelic sequence of $C^{\infty}$ models. We also define a bounded sequence 
of $C^{\infty}$ models. 

\begin{Definition}
A sequence of $C^{\infty}$ models $\{(\calX_i, \overline{\calL}_i)\}_{i=0}^{\infty}$ 
of $(X, L)$ is called an {\em adelic sequence} 
(resp. {\em a bounded sequence}) if it satisfies the following 
properties: 
\begin{enumerate}
\item[(i)] $\overline{\calL}_i$ is nef for all $i$;  
\item[(ii)]
There is a non-empty Zariski open set $U \subseteq \Spec(O_K)$ 
such that 
$\rest{\calX_i}{U} = \rest{\calX_j}{U}$ 
(which we denote by $\calX_U$) and $\rest{\calL_i}{U} = \rest{\calL_j}{U}$ in 
$\Pic(\calX_U)\otimes\QQ$ for all $i$ and $j$; 
The open set $U$ is called a {\em common base}; 
\item[(iii)]
For each $i, j$, there are a projective arithmetic variety 
$\calX_{i,j} \overset{\pi_{i,j}}{\longrightarrow} \Spec(O_K)$, 
birational morphisms 
$\mu_{i,j}^i: \calX_{i,j}\to\calX_i$ and 
$\mu_{i,j}^j: \calX_{i,j}\to\calX_j$, and a nef 
$C^{\infty}$ hermitian $\QQ$-line bundle $\overline{D_{i,j}}$ on $\Spec(O_K)$ such 
that 
\[
\pi_{i,j}^*(\overline{D_{i,j}}^{\otimes -1})
\precsim_{\pi_{i,j}^{-1}(U)} 
(\mu_{i,j}^i)^*(\overline{\calL}_i) \otimes 
(\mu_{i,j}^j)^*(\overline{\calL}_j^{\otimes -1})
\precsim_{\pi_{i,j}^{-1}(U)} 
\pi_{i,j}^*(\overline{D_{i,j}})
\]
and that 
$\adeg(\overline{D_{i,j}}) \to 0$ 
as $i, j \to +\infty$ (resp. $\adeg(\overline{D_{i,j}})$ is 
bounded with respect to $i, j$). 
\end{enumerate}
\end{Definition}

Let $\{(\calX_i, \overline{\calL}_i)\}_{n=0}^{\infty}$ be 
an adelic sequence (resp. bounded sequence) of $C^{\infty}$ models of $(X, L)$. 
Let $g: Y \to X$ is a finite morphism of projective varieties over $K$. 
Suppose, for each $n$, we are given a morphism $g_i: \calY_i \to \calX_i$ 
of projective arithmetic varieties over $\Spec(O_K)$ 
that extends $g: Y \to X$ such that 
$\rest{\calY_i}{U} = \rest{\calY_j}{U}$ and 
$\rest{g_i}{U} = \rest{g_j}{U}$ for all $i$ and $j$. 
As mentioned in \cite[p.~41]{Mor} we can then see that 
$\{(\calY_i, g_i^*(\overline{\calL}_i))\}_{n=0}^{\infty}$ 
is an adelic sequence (resp. bounded sequence) of $C^{\infty}$ models 
of $(Y, g^*(L))$. 

One has the following theorem (\cite[Theorem~(1.4)]{Zh} 
and \cite[Proposition~4.1.1]{Mor}). 

\begin{Theorem}[\cite{Zh}, \cite{Mor}]
\label{thm:adelic:intersection}
Let $\{(\calX_i, \overline{\calL}_i^{(n)})\}_{i=0}^{\infty}$
be adelic sequences of $C^{\infty}$ models of $(X, L)$ for 
$n = 1, 2, \ldots, \dim X+1$. 
Then the arithmetic intersection number 
\[
\adeg\left(
\acherncl_1\left( \overline{\calL}_i^{(1)}\right)
\acherncl_1\left( \overline{\calL}_i^{(2)}\right)
\cdots 
\acherncl_1\left( \overline{\calL}_i^{(\dim X + 1)}\right)
\right)
\] 
converges as $i$ tends to $+\infty$. 
\end{Theorem}

\medskip
For $i = 1, 2, \ldots$, let $f_i: X \to X$ be a morphism over $K$ with 
$f_i^*(L) \simeq L^{\otimes d_{f_i}}$ for some integer 
$d_{f_i} \geq 2$. 
Since we assume $L$ is ample, every $f_i$ is finite. 
We say that a sequence $\bd{f} = (f_i)_{i=1}^{\infty}$ 
is {\em adelically bounded} if 
there are a $C^{\infty}$ model  $(\calX, \overline{\calL})$ of $(X, L)$ 
such that $\overline{\calL}$ is nef, 
and a non-empty Zariski open set $U \subseteq \Spec(O_K)$ with the following 
properties: 
\begin{enumerate}
\item[(i)]
$f_i: X \to X$ extend to $f_{iU}: \calX_U \to \calX_U$ for each $i$, 
where $\calX_U = \rest{\calX}{U}$; 
\item[(ii)]
Let $\calZ_i$ be the normalization of 
$\calX_U \overset{f_{iU}}{\longrightarrow} 
\calX_U \hookrightarrow \calX$; We write $\widetilde{f}_i$ 
for the induced morphism $\widetilde{f}_i: \calZ_i \to \calX$, and  
set $\overline{\calN}_i = \widetilde{f}_i^*
\left(\overline{\calL}\right)^{\otimes \frac{1}{d_{f_i}}}$; 
We also put $\calZ_0 = \calX$ and $\overline{\calN}_0 = \overline{\calL}$;  
Then $\{(\calZ_i, \overline{\calN}_i)\}_{i=0}^{\infty}$ is 
a bounded sequence with the common base $U$. 
\end{enumerate} 

For example, it is easy to see that an iteration by a finite number 
of morphisms $g_1, \ldots, g_k$ in Example~\ref{eg:1} gives an 
adelically bounded sequence.  

Let us make an adelic sequence of $C^{\infty}$ models from 
an adelically bounded sequence. 

\begin{Proposition}
\label{prop:adelic:sequence}
Let $\bd{f} = (f_i)_{i=1}^{\infty}$ be adelically bounded, and 
$(\calX, \overline{\calL})$ and $U$ in the definition 
of an adelically bounded sequence.  
Let $\calX_i$ be the normalization of 
\[
\calX_U \overset{f_{1U}}{\longrightarrow} 
\calX_U \overset{f_{2U}}{\longrightarrow} \cdots 
\overset{f_{iU}}{\longrightarrow} 
\calX_U \hookrightarrow \calX.
\]
We write $F_i$ for the induced morphism $F_i: 
\calX_i \to \calX$, and  set $\overline{\calL}_i = 
F_i^*\left(\overline{\calL}\right)^{\otimes 
\frac{1}{d_{f_1}\cdots d_{f_i}}}$.  
We also put $\calX_0 = \calX$ and $\overline{\calL}_0 = \overline{\calL}$. 
Then $\{(\calX_i, \overline{\calL}_i)\}_{i=0}^{\infty}$ is 
an adelic sequence. 
\end{Proposition}

\Proof
Let $\{(\calZ_i, \overline{\calN}_i)\}_{i=0}^{\infty}$ 
be the bounded sequence of $C^{\infty}$ models and 
$\widetilde{f}_{i}: \calZ_i \to \calX$ the morphisms 
in the definition of an adelically bounded sequence. 

{\bf Step.~1}:\quad 
We compare $(\calX_i, \overline{\calL}_i)$ 
with $(\calX_{i+1}, \overline{\calL}_{i+1})$.
Since $\calZ_{i+1}$ is the normalization of 
$\calX_U \overset{f_{i+1 U}}{\longrightarrow} \calX_U 
\hookrightarrow \calX$, we see that 
$\calX_{i+1}$ is the normalization of 
\[
\calX_U \overset{f_{1U}}{\longrightarrow} 
\calX_U \overset{f_{2U}}{\longrightarrow} \cdots 
\overset{f_{iU}}{\longrightarrow} 
\calX_U \hookrightarrow \calZ_{i+1}. 
\]
We write $G_i: \calX_{i+1} \to \calZ_{i+1}$ for 
the induced morphism. Since $F_{i+1}$ and 
$\widetilde{f}_{i+1} \circ G_i$ are both morphisms form 
$\calX_{i+1}$ to $\calX$ and coincides over $U$, 
we get $F_{i+1} = \widetilde{f}_{i+1} \circ G_i$. 
Then 
\begin{equation}
\label{eqn:adelic:sequence}
\overline{\calL}_{i+1} 
= 
F_{i+1}^*\left(\overline{\calL}\right)^{\otimes 
\frac{1}{d_{f_1}\cdots d_{f_i} d_{f_{i+1}}}}
= 
(G_i)^*(\widetilde{f}_{i+1})^*\left(
\overline{\calL}\right)^{\otimes 
\frac{1}{d_{f_1}\cdots d_{f_i} d_{f_{i+1}}}}
= 
(G_i)^*\left(\overline{\calN_{i+1}}\right)^{\otimes 
\frac{1}{d_{f_1}\cdots d_{f_i}}}. 
\end{equation} 

Since $\{(\calZ_i, \overline{\calN}_i)\}_{i=0}^{\infty}$ is 
a bounded sequence, there are a projective arithmetic variety 
$\calW_{i+1} \overset{\pi_{i+1}}{\longrightarrow} \Spec(O_K)$, 
birational morphisms 
$\mu_{i+1}: \calW_{i+1}\to\calZ_0 = \calX$ and 
$\nu_{i+1}: \calW_{i+1}\to\calZ_{i+1}$, and a nef 
$C^{\infty}$ hermitian $\QQ$-line bundle $\overline{D_{i+1}}$ 
on $\Spec(O_K)$ such 
that 
\begin{equation}
\label{eqn:adelic:sequence:2}
\pi_{i+1}^*(\overline{D_{i+1}}^{\otimes -1})
\precsim_{\pi_{i+1}^{-1}(U)} 
\mu_{i+1}^*(\overline{\calN}_0) \otimes 
\nu_{i+1}^*(\overline{\calN}_{i+1}^{\otimes -1})
\precsim_{\pi_{i+1}^{-1}(U)} 
\pi_{i+1}^*(\overline{D_{i+1}})
\end{equation}
and that there is a constant $C$ such that 
$0 \leq  \adeg(\overline{D_{i+1}}) \leq C$ for all $i$.

There is a projective arithmetic variety 
$\calW_{i+1}^{\prime} \overset{\pi_{i+1}^{\prime}}{\longrightarrow} 
\Spec(O_K)$, 
birational morphisms 
$\mu_{i+1}^{\prime}: \calW_{i+1}^{\prime}\to\calX_{i}$ and 
$\nu_{i+1}^{\prime}: \calW_{i+1}^{\prime}\to\calX_{i+1}$ 
that are the identity map over $U$, and a morphism 
$H_i: \calW_{i+1}^{\prime} \to \calW_{i+1}$ 
that extends $f_{1U}\circ\cdots\circ f_{iU}$. 
Indeed, let $\calW_{i+1}^{\prime\prime}$ be the normalization of 
$\calX_U \overset{f_{1U}}{\longrightarrow} 
\calX_U \overset{f_{2U}}{\longrightarrow} \cdots 
\overset{f_{iU}}{\longrightarrow} 
\calX_U \hookrightarrow \calW_{i+1}$, and 
we may take $\calW_{i+1}^{\prime}$ as the Zariski closure 
of $\calX_U$ diagonally embedded in $\calX_i\times_{\Spec(O_K)}
\calW_{i+1}^{\prime\prime}\times_{\Spec(O_K)}\calX_{i+1}$. 
Thus we get the following diagram: 
\begin{equation*}
\begin{CD}
\calX_i @>{F_i}>> \calX_0 = \calZ_0 = \calX \\
@A{\mu_{i+1}^{\prime}}AA @A{\mu_{i+1}}AA  \\
\calW_{i+1}^{\prime} @>{H_i}>> \calW_{i+1} \\ 
@V{\nu_{i+1}^{\prime}}VV @V{\nu_{i+1}}VV   \\
\calX_{i+1} @>{G_i}>> \calZ_{i+1}.
\end{CD}
\end{equation*}
Since $F_i\circ \mu_{i+1}^{\prime}$ 
and $\mu_{i+1} \circ H_i$ are both morphisms 
from $\calW_{i+1}^{\prime}$ to $\calX_0$ and coincide over $U$, we get 
$F_i\circ \mu_{i+1}^{\prime}= \mu_{i+1} \circ H_i$. Similarly 
$G_i\circ \nu_{i+1}^{\prime}= \nu_{i+1} \circ H_i$. 
By \eqref{eqn:adelic:sequence}, we have
\begin{align*}
(\mu_{i+1}^{\prime})^* (\overline{\calL}_{i}) 
\otimes (\nu_{i+1}^{\prime})^*(\overline{\calL}_{i+1}^{\otimes -1}) 
& = 
(\mu_{i+1}^{\prime})^* \left(F_i^* \left(\overline{\calL_0}\right)^{\otimes 
\frac{1}{d_{f_1}\cdots d_{f_i}}}\right)
\otimes (\nu_{i+1}^{\prime})^*\left(G_i^* 
\left(\overline{\calN_{i+1}}\right)^{\otimes 
- \frac{1}{d_{f_1}\cdots d_{f_i}}}\right) \\
& = 
H_i^*\left(\mu_{i+1}^*\left(\overline{\calN}_0\right)^{\otimes 
\frac{1}{d_{f_1}\cdots d_{f_i}}}
 \otimes 
\nu_{i+1}^*\left(\overline{\calN}_{i+1}\right)^{\otimes - 
\frac{1}{d_{f_1}\cdots d_{f_i}}}\right). 
\end{align*}
Then by \eqref{eqn:adelic:sequence:2} we get 
\begin{equation}
\label{eqn:adelic:sequence:3}
(\pi_{i+1}^{\prime})^*(\overline{D_{i+1}})^{\otimes - 
\frac{1}{d_{f_1}\cdots d_{f_i}}}
\precsim_{\pi_{i+1}^{\prime -1}(U)} 
(\mu_{i+1}^{\prime})^* (\overline{\calL}_{i}) 
\otimes (\nu_{i+1}^{\prime})^*(\overline{\calL}_{i+1}^{\otimes -1}) 
\precsim_{\pi_{i+1}^{\prime -1}(U)} 
(\pi_{i+1}^{\prime})^*(\overline{D_{i+1}})^{\otimes 
\frac{1}{d_{f_1}\cdots d_{f_i}}}. 
\end{equation}

{\bf Step.~2}:\quad
For $i > j \geq 0$, we compare $(\calX_i, \overline{\calL}_i)$ 
with $(\calX_j, \overline{\calL}_{j})$.
We take a projective arithmetic variety 
$\calW \overset{\pi}{\longrightarrow} \Spec(O_K)$ 
such that, for each $\alpha = i+1, \ldots, j$,  
there is a birational morphism $p_{\alpha}: \calW\to\calW_{\alpha}^{\prime}$ 
that is the identity map over $U$. Indeed, we may take $\calW$ as the Zariski closure of $\calX_U$ diagonally embedded in 
$\calW_{i+1}^{\prime}\times_{\Spec(O_K)}
\calW_{i+2}^{\prime}\times_{\Spec(O_K)}\cdots \times_{\Spec(O_K)}
\calW_{j}^{\prime}$. For $\alpha = i+1, \ldots, j$, 
we define the birational morphism 
$\xi_{\alpha}: \calW\to\calX_{\alpha}$ by 
\begin{multline*}
\xi_i = \mu_{i+1}^{\prime}\circ p_{i+1}, \quad
\xi_{i+1} = \mu_{i+2}^{\prime}\circ p_{i+2} 
= \nu_{i+1}^{\prime}\circ p_{i+1}, \quad
\ldots, \\ 
\xi_{j-1} = \mu_{j}^{\prime}\circ p_{j} 
= \nu_{j-1}^{\prime}\circ p_{j-1}, \quad
\xi_{j} = \nu_{j}^{\prime}\circ p_{j}. 
\end{multline*}
Here $\mu_{\alpha+1}^{\prime}\circ p_{\alpha+1} 
= \nu_{\alpha}^{\prime}\circ p_{\alpha}$ 
for $\alpha = i+1, \ldots, j-1$
because $\mu_{\alpha+1}^{\prime}\circ p_{\alpha+1}$ and 
$\nu_{\alpha}^{\prime}\circ p_{\alpha}$ are the identity map 
over $U$ whence over $\Spec(O_K)$. 

In what follows, we write the group structure of the Picard group
additively. By \eqref{eqn:adelic:sequence:3}, we have 
\[
- \frac{1}{d_{f_1}\cdots d_{f_{\alpha}}} \pi^*(\overline{D_{\alpha+1}})
\precsim_{\pi^{-1}(U)} 
\xi_{\alpha}^* (\overline{\calL}_{\alpha}) 
- \xi_{\alpha+1}^*(\overline{\calL}_{\alpha+1}) 
\precsim_{\pi^{-1}(U)} 
\frac{1}{d_{f_1}\cdots d_{f_{\alpha}}} \pi^*(\overline{D_{\alpha+1}}). 
\]
Since $\xi_{i}^* (\overline{\calL}_{i}) 
- \xi_{j}^*(\overline{\calL}_{j}) 
= \sum_{\alpha=i}^{j-1} \left(
\xi_{\alpha}^* (\overline{\calL}_{\alpha}) 
- \xi_{\alpha+1}^*(\overline{\calL}_{\alpha+1}) \right)$, we get  
\[
- \sum_{\alpha=i}^{j-1} \frac{1}{d_{f_1}\cdots d_{f_{\alpha}}} \pi^*(\overline{D_{\alpha+1}})
\precsim_{\pi^{-1}(U)} 
\xi_{i}^* (\overline{\calL}_{i}) 
- \xi_{j}^*(\overline{\calL}_{j}) 
\precsim_{\pi^{-1}(U)} 
\sum_{\alpha=i}^{j-1} \frac{1}{d_{f_1}\cdots d_{f_{\alpha}}} \pi^*(\overline{D_{\alpha+1}}). 
\]
Since $0 \leq \adeg(\overline{D_{\alpha+1}}) \leq C$ for all $\alpha$, we have 
$\sum_{\alpha=i}^{j-1} \frac{1}{d_{f_1}\cdots d_{f_{\alpha}}} 
\adeg(\overline{D_{\alpha+1}}) \to 0$
as $i, j$ tends to $+\infty$. 
Thus $\{(\calX_i, \overline{\calL}_i)\}_{i=0}^{\infty}$ is an adelic sequence.
\QED

To give another construction of canonical heights 
for subvarieties, we need the following proposition 
(\cite[Proposition~2.1 and its proof]{Mor}). 

\begin{Proposition}[\cite{Mor}]
\label{prop:Mo:Prop:2:1}
Let $(\calX_1, \overline{\calL}_1)$ and $(\calX_2, \overline{\calL}_2)$ 
be two $C^{\infty}$ models of $(X, L)$ such that 
$\overline{\calL}_1$ and $\overline{\calL}_2$ are nef. 
Assume that there are 
a non-empty Zariski open set $U \subseteq \Spec(O_K)$, 
a projective arithmetic variety 
$\calX_3 \overset{\pi}{\longrightarrow} \Spec(O_K)$, 
birational morphisms $\mu_1: \calX_3\to\calX_1$ 
and $\mu_2: \calX_3\to\calX_2$, 
and a nef $C^{\infty}$ hermitian $\QQ$-line bundle 
$\overline{D}$ on $\Spec(O_K)$, such that 
\[
\pi^*(\overline{D}^{\otimes -1})
\precsim_{\pi^{-1}(U)} 
\mu_1^*(\overline{\calL}_1) \otimes 
\mu_2^*(\overline{\calL}_2^{\otimes -1})
\precsim_{\pi^{-1}(U)} 
\pi^*(\overline{D}). 
\]
Let $Y$ be a subvariety of $X_{\overline{K}}$, which is defined 
over a finite extension field $K'$ of $K$. Let $\calY$ be 
the Zariski closure of $\calX_3\times_{\Spec(O_K)} \Spec(O_K')$. 
Let $p: \calX_3\times_{\Spec(O_K)} \Spec(O_K') \to \calX_3$ be 
the natural morphism. Then 
\begin{multline*}
\left\vert
\adeg\left(\acherncl_1\left(\rest{(\mu_1\circ p)^*(\overline{\calL}_1)}{\calY}
\right)^{\cdot \dim Y +1}\right) - 
\adeg\left(\acherncl_1\left((\mu_2\circ p)^*(\rest{\overline{\calL}_2)}{\calY}
\right)^{\cdot \dim Y +1}\right)
\right\vert \\
\leq [K' :K] (\dim Y + 1) \deg(\rest{L}{Y}^{\dim Y}) 
\adeg(\overline{D}). 
\end{multline*}
\end{Proposition}

\begin{Proposition}
\label{prop:can:height:2}
Let $\bd{f} = (f_i)_{i=1}^{\infty}$ be adelically bounded, 
and $\{(\calX_i, \overline{\calL}_i)\}_{i=0}^{\infty}$ 
the adelic sequence of $C^{\infty}$ models in 
Proposition~\ref{prop:adelic:sequence}.  
Then for any subvariety $Y$ of $X_{\overline{K}}$, we have 
\begin{equation}
\label{eqn:can:height:2}
\widehat{h}_{L, \bd{f}}(Y) = 
\lim_{i\to\infty} h_{(\calX_i, \overline{\calL}_i)}(Y). 
\end{equation}
Moreover, the convergence is uniform with respect to $Y$. 
\end{Proposition}

\Proof
We first check that $\bd{f}$ belongs to $\widetilde{\calB}$ 
so that $\widehat{h}_{L, \bd{f}}(Y)$ is well-defined. 
Let $\{(\calZ_i, \overline{\calN}_i)\}_{i=0}^{\infty}$ 
be the bounded sequence of $C^{\infty}$ models in the definition of 
an adelically bounded sequence. 
By the projection formula, we have 
$h_{(\calZ_i, \overline{\calN}_i)}(Y)
= \frac{1}{d_{f_i}} h_{(\calZ_0, \overline{\calN}_0)}(f_i(Y))$. 
On the other hand, it follows from the projection formula, 
Proposition~\ref{prop:Mo:Prop:2:1} and the definition 
of an bounded sequence of $C^{\infty}$ models that 
there is a constant $C$ independent of $i$ and $Y$ such that 
\[
\left\vert
h_{(\calZ_i, \overline{\calN}_i)}(Y) 
- h_{(\calZ_0, \overline{\calN}_0)}(Y) 
\right\vert \leq C. 
\]
Since $(\calZ_0, \overline{\calN}_0) = (\calX, \overline{\calL})$, 
we get $\left\vert 
\frac{1}{d_{f_i}} h_{(\calX, \overline{\calL})}(f_i(Y)) 
- h_{(\calX, \overline{\calL})}(Y) 
\right\vert \leq C$ for all $i$ and $Y$. 
Thus $\bd{f}$ belongs to $\widetilde{\calB}$.
For the right-hand side of \eqref{eqn:can:height:2}, 
it follows from the projection formula that 
\[
h_{(\calX_i, \overline{\calL}_i)}(Y) = 
\frac{1}{d_{f_1}\cdots d_{f_i}}
h_{(\calX, \overline{\calL})}(f_i\circ\cdots\circ f_1(Y)). 
\]
Then by \eqref{eqn:can:height}, we get the first assertion.
Since 
\begin{multline*}
\left\vert h_{(\calX_i, \overline{\calL}_i)}(Y) - 
h_{(\calX_{i+1}, \overline{\calL}_{i+1})}(Y) \right\vert \\ 
= 
\left\vert 
\frac{1}{d_{f_1}\cdots d_{f_i}}
h_{(\calX, \overline{\calL})}(f_i\circ\cdots\circ f_1(Y))
- 
\frac{1}{d_{f_1}\cdots d_{f_{i+1}}}
h_{(\calX, \overline{\calL})}(f_{i+1}\circ\cdots\circ f_1(Y))
\right\vert \leq 
\frac{C}{d_{f_1}\cdots d_{f_i}},   
\end{multline*}
the convergence is uniform with respect to $Y$. 
\QED

\setcounter{equation}{0}
\section{Admissible metrics}
\label{sec:admissible:metrics}
In this section, we consider a local theory  
in the setting of \cite[\S2]{Zh}, and in case $X = \PP^N$ 
over $\CC$, we see its relation with Green currents in \cite{FW}.    

Let $K_{v}$ be an algebraically closed valuation field, 
$X$ a projective variety over $K_v$, and $L$ a line bundle on $X$. 
Let $\Vert\cdot\Vert$  be a continuous and bounded metric on $L$. 
(When $K_v$ is nonarchimedean, we refer to \cite[(1.1)]{Zh} for 
the definition of $\Vert\cdot\Vert$ being a continuous and bounded 
metric.)

We set 
\begin{equation*}
\overline{\calH} 
:= 
\left\{
\overline{f} := (f, \varphi) \;\left\vert
\begin{gathered}
\text{$f: X \to X$ is a morphism over $K_v$,} \\
\text{$\varphi: L^{\otimes d_{f}} \overset{\sim}{\rightarrow} f^*L$ 
with some integer $d_f \geq 2$}
\end{gathered}
\right.
\right\}.
\end{equation*}

For $\overline{f} = (f, \varphi) \in \overline{\calH}$, we set
\[
\overline{c}(\overline{f}) 
= \sup_{x \in X(K_v)} 
\left\vert
\log 
\frac{(\varphi^* f^*\Vert\cdot\Vert)^{\frac{1}{d_f}}}%
{\Vert\cdot\Vert} (x) 
\right\vert  \quad \in \RR. 
\]

Moreover, for $\bd{\overline{f}} 
= (\overline{f_i})_{i=1}^{\infty} \in \prod_{i=1}^{\infty} 
\overline{\calH}$, we set 
\[
\overline{c}(\bd{\overline{f}}) 
= \sup_{i \geq 1} \overline{c}(\overline{f_i}) 
\quad \in \RR \cup \{+\infty\}. 
\]

We define 
$\overline{\calB} 
= \left\{
\bd{\overline{f}} \in \prod_{i=1}^{\infty} \overline{\calH}
\;\left\vert\; \overline{c}(\bd{\overline{f}}) < +\infty \right.\right\}$. 
Note that the property that $\bd{\overline{f}}$ belongs 
to $\overline{\calB}$ is independent of the choice of 
bounded and continuous metrics $\Vert\cdot\Vert$. 

Let $S: \prod_{i=1}^{\infty} \overline{\calH} \to 
\prod_{i=1}^{\infty} \overline{\calH}$ be the shift map 
defined by $S((\overline{f_i})_{i=1}^{\infty}) = 
(\overline{f_{i+1}})_{i=1}^{\infty}$. 
It follows from $\overline{c}(S(\bd{\overline{f}})) 
\leq \overline{c}(\bd{\overline{f}})$ that
$S$ sends $\overline{\calB}$ to $\overline{\calB}$. 

\begin{Theorem}
\label{thm:admissible:metric}
Let $X$ be a projective variety over 
an algebraically closed valuation field $K_v$, 
and $L$ a line bundle on $X$. Let  
$\Vert\cdot\Vert$ be a continuous and bounded metric on $L$. 
Then there is a unique way to attach to each 
$\bd{\overline{f}} \in \overline{\calB}$ 
a metric $\canmet$ 
with the following properties\textup{:}
\begin{enumerate}
\item[(i)] 
The metric $\canmet$ 
is bounded and continuous such that  
\[
\sup_{x \in X(K_v)} 
\left\vert
\log \frac{\canmet}{\Vert\cdot\Vert} (x) 
\right\vert \leq 2 \overline{c}(\bd{\overline{f}})\textup{;}
\] 
\item[(ii)] 
$\varphi_1^* f_1^* \canmetS
= \canmet^{d_{f_1}}$. 
\end{enumerate}
Moreover, the metric $\canmet$ is independent of the choice 
of bounded and continuous metrics on $L$. 
We call $\canmet$ the admissible metric for 
$\bd{\overline{f}}$. 
\end{Theorem}

\Proof
We denote by $g_i$ the bounded and continuous function 
\[
g_i := 
\log\frac{\left(\varphi_i^*f_i^*\Vert\cdot\Vert\right)^{\frac{1}{d_{f_i}}}}%
{\Vert\cdot\Vert}
\]
on $X(K_v)$. We define the bounded and continuous metric 
$\Vert\cdot\Vert_i$ on $L$ by 
\begin{align*}
\Vert\cdot\Vert_0 & := \Vert\cdot\Vert, \\
\Vert\cdot\Vert_i & := 
\left(
\varphi_1^*f_1^* 
\cdots \left(\varphi_{i-1}^*f_{i-1}^* \left(
\varphi_{i}^*f_{i}^* \Vert\cdot\Vert
\right)^{\frac{1}{d_{f_i}}}\right)^{\frac{1}{d_{f_{i-1}}}}\cdots
\right)^{\frac{1}{d_{f_1}}} \qquad (\text{for $i \geq 1$}). 
\end{align*}
Then 
\[
\log \frac{\Vert\cdot\Vert_i}{\Vert\cdot\Vert_{i-1}} 
= 
\prod_{\beta=1}^{i-1} 
\left(\frac{1}{d_{f_{\beta}}}\varphi_{\beta}^*f_{\beta}^* \right)
\left[
\frac{1}{d_{f_{i}}}\varphi_{i}^*f_{i}^* (\log\Vert\cdot\Vert_0) 
- \log\Vert\cdot\Vert_0\right] 
= 
\prod_{\beta=1}^{i-1} 
\left(\frac{1}{d_{f_{\beta}}}\varphi_{\beta}^*f_{\beta}^* \right) 
g_i. 
\]
We then have 
\[
\log \frac{\Vert\cdot\Vert_i}{\Vert\cdot\Vert_{0}}
= \sum_{\alpha=1}^{i}  
\log \frac{\Vert\cdot\Vert_{\alpha}}{\Vert\cdot\Vert_{\alpha-1}}
= \sum_{\alpha=1}^{i}  
\prod_{\beta=1}^{\alpha-1} 
\left(\frac{1}{d_{f_{\beta}}}\varphi_{\beta}^*f_{\beta}^* \right) 
g_{\alpha}. 
\]
Since 
$\sup_{x \in X(K_v)} \vert g_{\alpha} (x)\vert \leq 
\overline{c}(\overline{\bd{f}})$ for any $\alpha \geq 1$, 
it follows that 
\begin{equation}
\label{eqn:admissible:metric}
\sup_{x \in X(K_v)} \left\vert 
\sum_{\alpha=1}^{i}  
\prod_{\beta=1}^{\alpha-1} 
\left(\frac{1}{d_{f_{\beta}}}\varphi_{\beta}^*f_{\beta}^* \right) 
g_{\alpha}(x) \right\vert
\leq 
\sum_{\alpha=1}^{i} 
\frac{\overline{c}(\overline{\bd{f}})}%
{\prod_{\beta=1}^{\alpha-1}d_{f_{\beta}}}
\leq 
\sum_{\alpha=1}^{i} \frac{\overline{c}(\overline{\bd{f}})}{2^{\alpha-1}}. 
\end{equation}
Thus, as $i$ tends to the infinity,  
$\sum_{\alpha=1}^{i}  
\prod_{\beta=1}^{\alpha-1} 
\left(\frac{1}{d_{f_{\beta}}}\varphi_{\beta}^*f_{\beta}^* \right) 
g_{\alpha}$ converges to  
a bounded and continuous function, which we denote by 
$g_{\bd{\overline{f}}}$. 
We set $\canmet := \Vert\cdot\Vert_0 \exp(g_{\bd{\overline{f}}})$.  

Let us check $\canmet$ satisfies (i) and (ii).  
Using \eqref{eqn:admissible:metric}, one finds 
\[
\sup_{x \in X(K_v)} 
\left\vert
\log \frac{\canmet}{\Vert\cdot\Vert_0} (x) 
\right\vert 
\leq 
\sum_{\alpha=1}^{\infty} \frac{\overline{c}(\bd{\overline{f}})}{2^{\alpha-1}}
\leq 2 \overline{c}(\bd{\overline{f}}). 
\]
This shows (i). To see (ii), we first note 
\[
\log \frac{\canmetS}{\Vert\cdot\Vert_0}
= 
\sum_{\alpha=1}^{\infty}  
\prod_{\beta=1}^{\alpha-1} 
\left(\frac{1}{d_{f_{\beta+1}}}\varphi_{\beta+1}^*f_{\beta+1}^* \right) 
g_{\alpha+1}. 
\]
Then we get 
\begin{align*}
\varphi_1^* f_1^* (\log \canmetS) 
& =
\varphi_1^* f_1^* (\log \Vert\cdot\Vert_0) 
+ \varphi_1^* f_1^* 
\left(\sum_{\alpha=1}^{\infty}  
\prod_{\beta=1}^{\alpha-1} 
\left(\frac{1}{d_{f_{\beta+1}}}\varphi_{\beta+1}^*f_{\beta+1}^* \right) 
g_{\alpha+1}\right) \\
& = 
\varphi_1^* f_1^* (\log \Vert\cdot\Vert_0) 
- d_{f_1} g_1 + 
d_{f_1} 
\left(\sum_{\alpha=1}^{\infty}  
\prod_{\beta=1}^{\alpha-1} 
\left(\frac{1}{d_{f_{\beta}}}\varphi_{\beta}^*f_{\beta}^* \right) 
g_{\alpha}\right) \\
& = 
\varphi_1^* f_1^* (\log \Vert\cdot\Vert_0) 
- d_{f_1} g_1 + d_{f_1} \left(\log \canmet - \log \Vert\cdot\Vert_0\right)
= d_{f_1} \log \canmet, 
\end{align*}
where we used 
$g_1 = \log\frac{\left(
\varphi_1^*f_1^*\Vert\cdot\Vert\right)^{\frac{1}{d_{f_1}}}}%
{\Vert\cdot\Vert}$ 
in the last equality. This shows (ii). 

Next we show the uniqueness. Suppose $\{\canmet^{\prime}\}$ 
are metrics on $L$ satisfying (i) and (ii). Then we have by (ii) 
\begin{multline*}
\sup_{x \in X(K_v)} 
\left\vert \log \frac{\canmet}{\canmet^{\prime}} (x) \right\vert 
=  
\sup_{x \in X(K_v)} 
\left\vert 
\left(\frac{1}{d_{f_1}} \varphi_1^*f_1^* \right)
\log \frac{\canmetS}{\canmetS^{\prime}} (x) \right\vert \\
= \cdots = 
\sup_{x \in X(K_v)} 
\left\vert 
\left(\frac{1}{d_{f_1}} \varphi_1^*f_1^* \right)\cdots
\left(\frac{1}{d_{f_i}} \varphi_i^*f_i^* \right)
\log \frac{\canmetSi}{\canmetSi^{\prime}} (x) \right\vert 
\leq 
\frac{1}{\prod_{\alpha=1}^{i} d_{f_{\alpha}}}
4 \overline{c}(S^i(\bd{\overline{f}})). 
\end{multline*}
Since $\overline{c}(S^i(\bd{\overline{{f}}})) \leq 
\overline{c}(\bd{\overline{f}})$, letting $i$ to the infinity, we find 
$\canmet = \canmet^{\prime}$. 

In the rest of this section, 
we show that $\canmet$ is independent of the choice of 
$\Vert\cdot\Vert$. Since the property of $\overline{f}\in\overline{\calB}$ 
is independent of the choice of $\Vert\cdot\Vert$, let us start 
from another bounded and continuous metric 
$\Vert\cdot\Vert^{\prime\prime}$ to obtain $\canmet^{\prime\prime}$. 
Set $\overline{c} = \sup_{x\in X(K_v)} \left\vert
\log\frac{\Vert\cdot\Vert^{\prime\prime}}{\Vert\cdot\Vert}
\right\vert$. Then 
\begin{align*}
\sup_{x\in X(K_v)} \left\vert
\log\frac{\canmetSi^{\prime\prime}}{\Vert\cdot\Vert}(x)
\right\vert
& \leq 
\sup_{x\in X(K_v)} \left\vert
\log\frac{\canmetSi^{\prime\prime}}{\Vert\cdot\Vert^{\prime\prime}}(x)
\right\vert + \overline{c} \\
& \leq 
2 \overline{c}(S^i(\overline{\bd{f}}))
+ 5 \overline{c}
\leq 
2 \overline{c}(\overline{\bd{f}}) + 5 \overline{c}
\end{align*}
for any $i \geq 1$. Now the above argument of uniqueness of 
$\canmet$ implies $\canmet^{\prime\prime} = \canmet$. 
This completes the proof. 
\QED

In the rest of this section, we consider a case $K_v=\CC$, $X=\PP^N$ and 
$L= \OO_{\PP^N}(1)$. We endow $\OO_{\PP^N}(1)$ 
with the Fubini-Study metric $\Vert\cdot\Vert_{FS}$. 

In this case, to give $\overline{f}=(f,\varphi)\in\overline{\calH}$ is
equivalent to give homogeneous polynomials $F_0(X_0,\cdots,X_N),
\ldots, F_N(X_0,\cdots,X_N) \in \CC[X_0,\cdots,X_N]$ of degree $d_f$
such that $f=(F_0:\cdots:F_N)$. Indeed $F_k=\varphi^*f^*(X_k)$, where
$X_k$ is regarded as an element of $\Gamma(\PP^N, \OO_{\PP^N}(1))$.
Note that, since $f: \PP^N\to\PP^N$ is a morphism, the only common
zero of $F_0, \ldots, F_N$ is $0\in\CC^{N+1}$.

Let $\bd{\overline{f}}=(\overline{f_i})_{i=1}^{\infty} \in
\overline{\calB}$. Let $F_0^i(X_0,\cdots,X_N), \ldots,
F_N^i(X_0,\cdots,X_N) \in \CC[X_0,\cdots,X_N]$ be homogeneous
polynomials of degree $d_{f_i}$ corresponding to $\overline{f_i}$. We
put, for $i\geq 1$,
\[
F^i := (F_0^i, \cdots, F_N^i) : 
\CC^{N+1}\to \CC^{N+1}.
\] 

For $x=(x_0, \cdots, x_N) \in \CC^{N+1}$, we put $|x| = \sqrt{|x_0|^2
  + \cdots + |x_N|^2}$. 

We set $G_0(x) = \log |x|^2$, and for $i \geq 1$ 
we define $G_i: \CC^{N+1}\setminus\{0\}\to\RR$ by
\[
G_i(x) = \frac{1}{\prod_{\alpha=1}^{i} d_{f_{\alpha}}} 
\log | F^i\circ\cdots\circ F^1(x) |^2.
\]

Let $\calF_{\bd{\overline{f}}} \subset \PP^N(\CC)$ be the largest open
set on which the family $\{f_1,\; f_2\circ f_1,\; f_3\circ f_2\circ f_1,\;
\cdots\}$ is normal, and $\calJ_{\bd{\overline{f}}}$ be its complement.

\begin{Lemma}
\label{lemma:estimate:Green:function}
For any $i \geq1$ and $x=(x_0, \cdots, x_N) \in
\CC^{N+1}\setminus\{0\}$, we have
\[
|G_i(x) - G_{i-1}(x)| \leq 
\frac{\overline{c}(\bd{\overline{f}})}{2^{i-2}}.
\]
\end{Lemma} 

\Proof 
As in the proof of Theorem~\ref{thm:admissible:metric}, 
we let $g_i = \log
\frac{(\varphi_i^*f_i^*
\Vert\cdot\Vert_{FS})^{\frac{1}{d_{f_i}}}}{\Vert\cdot\Vert_{FS}}$. 
For $x=(x_0, \cdots, x_N) \in \CC^{N+1}\setminus\{0\}$ and 
$\pi(x)=(x_0:\cdots:x_N) \in \PP^N(\CC)$, we 
auxiliary take $k$ such that $F_k^i(x)\neq 0$. We regard $F_k^i$ as an
element of $\Gamma(\PP^N, \OO_{\PP^N}(d_{f_i}))$. Noting that
$\varphi_i^*f_i^*(X_k) = F_k^i$, we then get
\[
g_i(\pi(x))^{d_{f_i}} = 
\log \frac{\frac{|F_k^i(x)|}{|F^i(x)|}}{\frac{|F_k^i(x)|}{|x|^{d_{f_i}}}}
= \log \frac{|x|^{d_{f_i}}}{|F^i(x)|}. 
\]
It follows from $|g_i(\pi(x))| \leq \overline{c}(\bd{\overline{f}})$ that 
\[
\left\vert \frac{1}{d_{f_i}} \log |F^i(x)| - \log |x| \right\vert
\leq \overline{c}(\bd{\overline{f}}).
\]
Then we have
\[
\left\vert \frac{1}{d_{f_i}}\log |F^i(F^{i-1}\circ\cdots\circ F^1(x))| 
- \log |F^{i-1}\circ\cdots\circ F^1(x)| \right\vert
\leq \overline{c}(\bd{\overline{f}}).
\]
Thus
\[
|G_i(x) - G_{i-1}(x)| \leq 
\frac{2 \overline{c}(\bd{\overline{f}})}%
{\prod_{\alpha=1}^{i-1} d_{f_{\alpha}}} \leq 
\frac{\overline{c}(\bd{\overline{f}})}{2^{i-2}}.
\]
\QED

\begin{Proposition}[cf. \cite{FW}]
\label{prop:estimate:Green:function}
Let $X =\PP^N$ over $\CC$, $(L, \Vert\cdot\Vert) = (\OO_{\PP^N}(1),
\Vert\cdot\Vert_{FS})$, and $\bd{\overline{f}}\in \overline{\calB}$.
\begin{enumerate}
\item[(1)] 
As $l$ tends to the infinity, $G_l$ converges uniformly to
a plurisubharmonic and continuous function $G_{\bd{\overline{f}}}:
\CC^{N+1}\setminus\{0\} \to\RR$.
\item[(2)]
Let $T_{\bd{\overline{f}}}$ be the positive closed $(1,1)$ current
on $\PP^N(\CC)$ such that $\pi^* T_{\bd{\overline{f}}} = d d^c
(G_{\bd{\overline{f}}})$, where $\pi: \CC^{N+1}\setminus\{0\}
\to\PP^N$ is the natural projection. Then
$\Supp(T_{\bd{\overline{f}}}) = \calJ_{\bd{\overline{f}}}$.
\end{enumerate}
We call $G_{\bd{\overline{f}}}$ the Green function for 
$\bd{\overline{f}}$, and $T_{\bd{\overline{f}}}$ the  
Green current for $\bd{\overline{f}}$. 
\end{Proposition}

\Proof 
Using Lemma~\ref{lemma:estimate:Green:function}, one can prove
(1) as in \cite[\S2]{FW}, and (2) as in \cite[Proposition~8]{FW},
\cite[Th\'eor\`eme~1.6.5]{Si}.  
\QED

Let $\canmetP$ be the admissible metric for $\overline{\bd{f}}$ 
on $\OO_{\PP^N}(1)$ over
$\PP^N_{\CC}$. We define the first Chern current of $(\OO_{\PP^N}(1),
\canmetP)$ by
\[
\cherncl_1(\OO_{\PP^N}(1), \canmetP) := 
d d^c [-\log \widehat{\Vert s \Vert}_{\OO_{\PP^N(1)}, \bd{\overline{f}}}^2] 
+ \delta_{\zero(s)}
\]
for any nonzero section $s$ of $\OO_{\PP^N}(1)$.  It is noted in
\cite[Proposition~3.3.1]{Ka} that, 
when $f_1=f_2=\cdots (=f)$, the first Chern current
of $\OO_{\PP^N}(1)$ with the admissible metric coincides with the
Green current. The following lemma says this also holds for
$\bd{\overline{f}}\in \overline{\calB}$.

\begin{Proposition}
\label{prop:Chern:current}
With the notation and assumption as in
Proposition~\ref{prop:estimate:Green:function}, we have 
\[
\cherncl_1(\OO_{\PP^N}(1), \canmetP) = T_{\bd{\overline{f}}}.
\]
\end{Proposition}

\Proof 
Let $U_0 := \{(1:x_1:\cdots:x_N)\}$ be an open set of
$\PP^N(\CC)$.  Since it suffices to prove the equality locally, we
show it over $U_0$. (The case over $U_i :=\{(x_0:x_1:\cdots:x_N) \mid
x_i\neq 0\}$ is proven similarly.) We take $(x_1,\cdots,x_N)$ as a
coordinate of $U_0$. Then for $x=(1,x_1,\cdots,x_N) \in \CC^{N+1}$
and $\pi(x)=(1:x_1:\cdots:x_N) \in \PP^N(\CC)$, we have
\begin{align*}
& \log\frac{\canmetP^2}{\Vert\cdot\Vert_{FS}^2}(\pi(x))
= 
2 \sum_{\alpha=1}^{\infty}
\prod_{\beta=1}^{\alpha-1}
\left(\frac{1}{d_{f_{\beta}}} \varphi_{\beta}^* f_{\beta}^* \right)
g_{\alpha}(\pi(x)) \\
& \quad = 
\sum_{\alpha=1}^{\infty} 
\frac{1}{\prod_{\beta=1}^{\alpha-1} d_{f_{\beta}}}
\left(
\log \left\vert F^{\alpha-1}\circ\cdots\circ F^1(x)
\right\vert^2
- 
\frac{1}{d_{f_{\alpha}}}
\log \left\vert F^{\alpha}\circ F^{\alpha-1}\circ\cdots\circ F^1(x)
\right\vert^2 
\right) \\
& \quad = 
\sum_{\alpha=1}^{\infty} 
\left( G_{\alpha-1}(x) - G_{\alpha}(x) \right)
= \log \vert x \vert^2 - 
G_{\bd{\overline{f}}}(x).
\end{align*}
Taking $d d^c$ on both sides, we get the assertion.  
\QED

We remark the following lemma. 

\begin{Lemma}
\label{lemma:no:overline}
Let $X =\PP^N$ over $\CC$, $(L, \Vert\cdot\Vert) = (\OO_{\PP^N}(1),
\Vert\cdot\Vert_{FS})$, and 
$\bd{\overline{f}} = \left((f_i, \varphi_i)\right)_{i=0}^{\infty}, 
\bd{\overline{f}}' = \left((f_i^{'}, \varphi_i^{'})\right)_{i=0}^{\infty} 
 \in \overline{\calB}$. If $f_i = f_i^{'}$ for all $i$, then 
 $T_{\bd{\overline{f}}} = T_{\bd{\overline{f}}'}$
\end{Lemma}

\Proof
Let 
$F^{i'} = (F_0^{i'}, \cdots, F_N^{i'}) : \CC^{N+1}\to \CC^{N+1}$ 
be the lift of 
$f_i^{'}: \PP^N \to \PP^N$ given by $\varphi_i^{'}$. 
Then $F^{i'} = F^{i}$ for some non-zero constant $c_i \in \CC$. Then 
\begin{align*}
G_i^{'}(x) 
& = \frac{1}{\prod_{\alpha=1}^{i} d_{f_{\alpha}}} 
\log | F^{i'}\circ\cdots\circ F^{1'}(x) |^2 \\
& = \frac{1}{\prod_{\alpha=1}^{i} d_{f_{\alpha}}} 
\log | F^i\circ\cdots\circ F^1(x) |^2 + 
\frac{1}{\prod_{\alpha=1}^{i} d_{f_{\alpha}}}  
\log \left\vert 
c_1^{d_{f_i}\cdots d_{f_2}} c_2^{d_{f_i}\cdots d_{f_3}} 
\cdots c_{i-1}^{d_{f_i}} 
\right\vert^2 \\
& = 
G_i(x) + \sum_{k=1}^{i-1} \frac{1}{\prod_{\beta=1}^{k} d_{f_{\beta}}} 
\log | c_k |^2.
\end{align*}
By Proposition~\ref{prop:estimate:Green:function}(1), we get 
$G_{\bd{\overline{f}}} = 
G_{\bd{\overline{f}}'} + 
\sum_{k=1}^{\infty} \frac{1}{\prod_{\beta=1}^{k} d_{f_{\beta}}} 
\log | c_k |^2$ (The last term converges).  
Then Proposition~\ref{prop:estimate:Green:function}(2) yields the lemma. 
\QED

\setcounter{equation}{0}
\section{Equidistribution of small points on $\PP^1$}
\label{sec:equidistribution}
In this section, we show equidistribution of small points on $\PP^1$ 
for sequences of morphisms. 
Let us first recall some facts on analytic torsions. 
For the details, we refer to \cite{SABK}. 

Let $(M, k)$ be a $d$-dimensional compact K\"ahler manifold with the 
K\"ahler metric $k$. Let $(L, h)$ a $C^{\infty}$ hermitian line bundle 
over $M$.  The vector space $A^{0,q}(L^{\otimes n})$ of smooth $(0,q)$ 
forms on $M$ with values in $L^{\otimes n}$ is equipped with the $L^2$-metric 
given by 
\[
(s, t)_{L^2} 
= \frac{1}{(2 \pi)^d} \int_M \langle s(x), t(x)\rangle dv 
\qquad (s, t \in A^{0,q}(L^{\otimes n})), 
\]
where $dv$ is the normalized volume form on $M$ associated 
with $k$. Let $\overline{\partial}: A^{0,q}(L^{\otimes n}) \to 
A^{0,q+1}(L^{\otimes n})$ be the Dolbeault operator and 
$\overline{\partial}^*: A^{0,q}(L^{\otimes n}) \to 
A^{0,q-1}(L^{\otimes n})$ its formal adjoint with respect to 
the $L^2$-metric. Let $\square_n^q := (\overline{\partial} + 
\overline{\partial}^*)^2$ be the Laplacian acting on 
$A^{0,q}(L^{\otimes n})$. 
Let $\sigma(\square_n^q)$ be the spectrum of  
$\square_n^q$. We set, for $s \in \CC$, 
\[
\zeta_{q, n}(s) := \sum_{\lambda \in 
\sigma(\square_n^q)\setminus\{0\}}\lambda^{-s}. 
\]
The zeta function $\zeta_{q, n}(s)$ converges absolutely for 
$\Rea s >d$, extends meromorphically to 
the whole plane, and is holomorphic at $s=0$. Then (the twice of the 
logarithm) of the Ray-Singer analytic torsion $T(L^{\otimes n})$ 
is defined by, e.g., \cite[p.~132]{SABK}, 
\[
T(L^{\otimes n}) = \sum_{q=0}^d (-1)^{q+1} q \zeta^{\prime}_{q, n}(0). 
\]

Vojta \cite[Propositon~2.7.6]{Vo} proved the following theorem 
(see also Bismut--Vasserot \cite{BV}). We remark that 
in \cite{Vo}, $\tau(E_{\sigma})$ is defined by 
$\tau(E_{\sigma}) = - T(E_{\sigma})$. 

\begin{Theorem}[\cite{Vo}]
\label{thm:Vojta}
For each $q$ \textup{(}$0 \leq q \leq d$\textup{)}, 
there exists a constant $C$ such that, 
as $n \to \infty$, 
\[
\zeta^{\prime}_{q, n}(0) \geq - C n^d \log n. 
\]
\end{Theorem}

If $d=1$ then $T(L^{\otimes n}) = \zeta^{\prime}_{1, n}(0)$. Hence, 
if $M$ is a compact Riemann surface, then 
\begin{equation}
\label{eqn:Vojta}
T(L^{\otimes n}) \geq - C n \log n.
\end{equation}

\begin{Question}
\label{quesiton:analytic:torsion}
Let $(M, k)$ be a $d$-dimensional compact K\"ahler manifold with the 
K\"ahler metric $k$, and $L$ a positive line bundle over $M$. 
Let $h$ be an arbitrary $C^{\infty}$ metric 
on $L$. (We do not assume $\cherncl_1(L,h)$ is positive). Then 
is it true that, for any positive $\varepsilon$, there is $n_0$ such 
that $T(L^{\otimes n}) \geq - \varepsilon n^{d+1}$ for all $n \geq n_0$?  
\end{Question}

It follows from \eqref{eqn:Vojta} that 
Question~\ref{quesiton:analytic:torsion} is true for $d=1$. 
For any $d \geq 1$, 
Bismut--Vasserot \cite{BV} showed that 
if $\cherncl_1(L,h)$ is positive then 
$T(L^{\otimes n}) = O(n^d \log d)$. 
They also showed that $T(L^{\otimes n}) = O(n^{d+1})$ for 
an arbitrary $C^{\infty}$ metric $h$ (\cite[Theorem~11]{BV}). 
Question~\ref{quesiton:analytic:torsion} seems interesting, 
because the validity of Question~\ref{quesiton:analytic:torsion} 
for $d=1$ is the key to prove equidistribution of small points 
on $\PP^1$ in Theorem~\ref{thm:equidistribution} below.  

\bigskip
Before going back to sequences of morphisms, 
we need some results on arithmetic surfaces. 

Let $K$ be a number field, and $O_K$ its ring of integers. 
Let $\calX$ be a projective arithmetic surface (i.e., a projective 
arithmetic variety of dimension $2$), and 
$\overline{\calL} = (\calL, \{\Vert\cdot\Vert_{\sigma}\}_{\sigma})$ 
a $C^{\infty}$ hermitian line bundle. 

We set $\Gamma = H^0(\calX, {\calL})$. 
Then $\Gamma$ is a lattice of 
$V := H^0(\calX, {\calL}) \otimes_{\ZZ} \RR$. 
By fixing an isomorphism $V \simeq \RR^{\oplus \rank \Gamma}$, 
we give the Lebesgue measure $\mu$ on $V$. 
Since $V$ is the complex-conjugation-invariant subspace of 
$H^0(\calX, \calL) \otimes_{\ZZ} \CC$, 
the metrics $\{\Vert\cdot\Vert_{\sigma}\}_{\sigma}$ induces 
the $L^2$-norm and $\sup$ norm on $V$. 
We set 
$\displaystyle{\vol_{L^2}(\Gamma) = \frac{\mu(V / \Gamma)}{\mu(B_{L^2}(V))}}$ 
and $\displaystyle{\vol_{\sup}(\Gamma) 
= \frac{\mu(V / \Gamma)}{\mu(B_{\sup}(V))}}$, 
where $B_{L^2}(V) = \{x \in V \mid \Vert x\Vert_{L^2} \leq 1\}$ 
and $B_{\sup}(V) = \{x \in V \mid \Vert x\Vert_{\sup} \leq 1\}$. 
The arithmetic Euler characteristics 
$\chi_{L^2}(\calX, \overline{\calL})$ and 
$\chi_{\sup}(\calX, \overline{\calL})$ are then 
respectively defined by 
\[
\chi_{L^2}(\calX, \overline{\calL}) 
= - \log \vol_{L^2}(\Gamma)
\quad\text{and}\quad 
\chi_{\sup}(\calX, \overline{\calL}) 
= - \log \vol_{\sup}(\Gamma). 
\]

\begin{Theorem}
\label{thm:small:section}
Let $\calX$ be a projective arithmetic surface over $\Spec(O_K)$,
and $\overline{\calL} = \left(\calL, 
\{\Vert\cdot\Vert_{\sigma}\}_{\sigma}\right)$ 
a $C^{\infty}$ hermitian line bundle on $\calX$. 
Assume $\calL$ is ample on $\calX$. 
\begin{enumerate}
\item[(1)] 
If $\adeg\left(\acherncl_1(\overline{\calL})^2\right) > 0$,  
then there is a constant $C$ such that, as $n \to \infty$,  
\[
\chi_{\sup}(\calX, \overline{\calL}^{\otimes n})
\geq  
\frac{\adeg\left(\acherncl_1(\overline{\calL})^2\right)}{2} n^2 
+ C n \log n. 
\]
\item[(2)]
For any $\varepsilon > 0$, there exist a positive 
integer $n$ and a non-zero section $s \in H^0(\calX, \calL^{\otimes n})$ 
such that 
\[
\Vert s \Vert_{\sigma, \sup} \leq 
\exp\left(n \left(\varepsilon - 
\frac{\adeg\left(\acherncl_1(\overline{\calL})^2\right)}%
{2 [K:\QQ] \deg(\calL\otimes_{O_K} K)} \right) \right)
\] 
for all embeddings $\sigma: K \hookrightarrow \CC$. 
\end{enumerate}
\end{Theorem}

\Proof
If $\cherncl_1(\overline{\calL} \otimes_{O_K^{\sigma}} \CC)$ is positive 
for all $\sigma$, then Theorem~\ref{thm:small:section} is just 
a very special case of \cite[Theorem~8, Theorem~9]{GS},  
\cite[Chap.~VIII]{SABK}. 
Using \eqref{eqn:Vojta}, one can prove Theorem~\ref{thm:small:section} 
as in {\it op.\! cit}. For the sake of completeness, 
we sketch a proof here. 

(1) 
Let $B_r$ denote the unit ball in $\RR^r$ with 
the Euclidean metric. We have, for large $n$,  
\begin{align*}
\chi_{\sup}(\calX, \overline{\calL}^{\otimes n}) \\
& = \chi_{L^2}(\calX, \overline{\calL}^{\otimes n}) + O(n \log n) \\
& = \adeg\left(\acherncl_1\left(\lambda(\overline{\calL}^{\otimes n}), 
\text{Quillen metric}\right)\right) -\sum_{q =1}^2 (-1)^q 
\log \# H^q(\calX, \calL^{\otimes n}) \\
& \qquad + \log \text{vol}(B_{\rank H^0(\calX, \calL^{\otimes n})}) 
+ \frac{1}{2} T(\calL^{\otimes n} \otimes_\ZZ \CC) + O(n \log n) \\
& = \frac{\adeg\left(\acherncl_1(\overline{\calL})^2\right)}{2} n^2 
+ \frac{1}{2}  T(\calL^{\otimes n} \otimes_\ZZ \CC) + O(n \log n) \\
& \geq 
\frac{\adeg\left(\acherncl_1(\overline{\calL})^2\right)}{2} n^2 
+ C n \log n. 
\end{align*}
Here we used the $L^2$-$\sup$ comparison due to Gromov 
(cf. \cite[Lemma~30]{GS}) in the first equality, the definition 
of Quillen metrics (cf. \cite[Chap.~VIII.~Lemma~1]{SABK}, \cite[\S 4]{GS}) 
in the second equality, Gillet-Soul\'e's arithmetic Riemann-Roch theorem 
(cf. \cite[Theorem~7]{GS}, \cite[Chap.~VIII]{SABK}) in the third equality, 
and Vojta's estimate \eqref{eqn:Vojta} in the last inequality.   

(2) Let $\varepsilon$ be any positive number. 
Let $\overline{\calL}^{\prime}$ be the $\QQ$-line bundle 
$\calL$ on $\calX$ equipped with the metric 
$\Vert\cdot\Vert_{\sigma}^{\prime}$ given by 
\begin{equation}
\label{eqn:metric:on:L:prime}
\Vert\cdot\Vert_{\sigma}^{\prime} = 
\Vert\cdot\Vert_{\sigma} \exp\left( 
\frac{\adeg\left(\acherncl_1(\overline{\calL})^2\right)}%
{2 [K:\QQ] \deg(\calL\otimes_{O_K} K)} - \varepsilon \right)
\end{equation}
for all $\sigma: K\hookrightarrow \CC$. 
Then 
\begin{align*}
\adeg\left(\acherncl_1(\overline{\calL^{\prime}})^2\right) 
& = 
\adeg\left(\acherncl_1(\overline{\calL})^2\right)
+ 2 [K:\QQ] \deg(\calL\otimes_{O_K} K) 
\left(\varepsilon - 
\frac{\adeg\left(\acherncl_1(\overline{\calL})^2\right)}%
{2 [K:\QQ] \deg(\calL\otimes_{O_K} K)}\right)\\
& = 
2 [K:\QQ] \deg(\calL\otimes_{O_K} K)  \varepsilon > 0. 
\end{align*}

We set $\Gamma_n = H^0(\calX, {\calL}^{\otimes n})$ 
and  $V_n = H^0(\calX, {\calL}^{\otimes n}) \otimes_{\ZZ} \RR$. 
We endow $V_n$ with the sup metric $\Vert\cdot\Vert_{\sup}^{\prime}$ 
induced by $\{\Vert\cdot\Vert_{\sigma}^{\prime n}\}_{\sigma}$. 
We set $\lambda_1(\Gamma_n) 
= \{ \Vert s \Vert_{\sup}^{\prime} \mid 0 \neq s \in \Gamma_n\}$. 
We also set 
$\vol_{\sup}(\Gamma_n) = \frac{\mu(V_n /\Gamma_n)}{\mu(B_{\sup}(V_n))}$, 
where $\mu$ is a Lebesgue measure on $V_n$ and 
$B_{\sup}(V_n) = \{x \in V_n \mid \Vert x \Vert_{\sup}^{\prime} \leq 1\}$. 
Minkowski's convex body theorem states that 
\[
(\dim V_n) \log \lambda_1(\Gamma_n) 
\leq (\dim V_n) \log 2 + \log \vol_{\sup}(\Gamma_n). 
\] 
It follows from (1) that $\log \vol_{\sup}(\Gamma_n) \leq 
- \frac{\adeg\left(\acherncl_1(\overline{\calL}^{\prime})^2\right)}{2} n^2 
- C n \log n$ for large $n$. 
On the other hand, $\dim V_n = [K:\QQ] \left\{\deg(\calL\otimes_{O_K} K) n 
+ \left(1- g(\calX_{\sigma}(\CC))\right)\right\}$ for large $n$. 
Hence, for any positive $\varepsilon^{\prime} > 0$, we find 
$\log \lambda_1(\Gamma_n) \leq n \left(\varepsilon^{\prime} - 
\frac{\adeg\left(\acherncl_1(\overline{\calL}^{\prime})^2\right)}%
{2 [K:\QQ] \deg(\calL\otimes_{O_K} K)}\right)$ for large $n$. 
Taking $\varepsilon^{\prime} = \varepsilon$, we see 
$\log \lambda_1(\Gamma_n) \leq 0$. Then 
the assertion follows from \eqref{eqn:metric:on:L:prime}. 
\QED

\bigskip
In what follows, we fix an embedding 
$\sigma_0: K \hookrightarrow \CC$, and extend $\sigma_0$ to 
$\sigma_0: \overline{K} \hookrightarrow \CC$. 

We consider the case $X = \PP^1$ and $L = \OO_{\PP^1}(1)$. 
For $i = 1, 2, \ldots$, 
let $f_i: \PP^1 \to \PP^1$ be a surjective morphism over $K$ of 
degree $d_{f_i} \geq 2$. We define the morphism 
$f_{i\CC}: \PP^1_{\CC} \to \PP^1_{\CC}$ over $\CC$ 
by $f_{i\CC} := f_i \times_{K^{\sigma_0}} \CC: 
\PP^1 \times_{K^{\sigma_0}} \CC \to 
\PP^1 \times_{K^{\sigma_0}} \CC$. 

For $\bd{f} := (f_i)_{i=1}^{\infty}$, we assume the following 
two conditions: 
\begin{enumerate}
\item[(C-i)]
$\bd{f}$ is adelically bounded such that $\calL_0$ is ample on 
$\calX_0$ (cf. 
\S\ref{sec:canonical:heights:for:subvarieties:2});  
\item[(C-ii)]
There are isomorphisms $\varphi_i: \OO_{\PP^1_{\CC}}(d_{f_i}) 
\to f_{i\CC}^{*} \OO_{\PP^1_{\CC}}(1)$ over $\CC$ for $i = 1, 2, \ldots$ 
such that 
\[
 \sup_{i \geq 1} \sup_{x \in \PP^1(\CC)} 
\left\vert
\log 
\frac{(\varphi_i^* f_{i\CC}^*\Vert\cdot\Vert_{FS})^{\frac{1}{d_{f_i}}}}%
{\Vert\cdot\Vert_{FS}} (x) 
\right\vert  < +\infty, 
\]
where $\Vert\cdot\Vert_{FS}$ is the Fubini-Study metric 
on $\OO_{\PP^1_{\CC}}(1)$ over $\PP^1_{\CC}$. 
\end{enumerate}

For example, iterations by a finite number of morphisms 
$g_1, \ldots, g_k$ in Example~\ref{eg:1} satisfy these 
conditions. For $\bd{f}$ satisfying (C-i) and (C-ii), 
it follows from Proposition~\ref{prop:can:height:2} and 
Proposition~\ref{prop:estimate:Green:function} that 
one has the canonical height function 
$\widehat{h}_{\PP^1, \bd{f}}: \PP^1(\overline{K}) \to \RR$ 
and the Green current $T_{\bd{f}}$ on $\PP^1(\CC)$. 
Here, in virtue of Lemma~\ref{lemma:no:overline}, we use the notation 
$T_{\bd{f}}$ instead of $T_{\overline{\bd{f}}}$, where 
$\overline{\bd{f}} = \left((f_i, \varphi_i)\right)_{i=1}^{\infty}$. 

Let $(x_j)_{j=1}^{\infty}$ be a sequence with $x_j \in \PP^1(\overline{K})$. 
The sequence $(x_j)_{j=1}^{\infty}$ is said to be {\em generic} 
if, for any $x \in \PP^1(\overline{K})$, there exists $j_0$ such that 
$x_j \neq x$ for any $j \geq j_0$. 
The sequence $(x_j)_{j=1}^{\infty}$ is called a sequence of {\em small} 
points with respect to $\widehat{h}_{\bd{f}}$ if 
$\lim_{j \to \infty} \widehat{h}_{\bd{f}}(x_j) = 0$. 

For $x \in \PP^1(\overline{K})$, let 
$G(x) := \{ \Gal(\overline{K}/K)\cdot x\} \subset \PP^1(\overline{K})$ 
be the orbit of $x$ under the Galois group.  
For $y \in \PP^1(\overline{K})$,  
let $\delta_y$ denote the Dirac measure of mass $1$ on $\PP^1(\CC)$ 
supported in $y$, where we regard $y$ as a point of $\PP^1(\CC)$ through 
the embedding $\sigma_0: \overline{K} \hookrightarrow \CC$. 

\begin{Example}
\label{eg:small:points}
We fix $a \in \PP^1(\overline{K})$ with 
$\widehat{h}_{\bd{f}}(a) \neq 0$. Let $x_j\in \PP^1(\overline{K})$ 
be a point with $f_j\circ f_{j-1}\circ \cdots \circ f_1(x_j)=a$.
Since 
$\widehat{h}_{\bd{f}}(x_j) = \frac{1}{d_{f_1}\cdots d_{f_j}} 
\widehat{h}_{\bd{f}}(a)$, $x_j$'s are distinct points such that 
$\lim_{j \to \infty} \widehat{h}_{\bd{f}}(x_j) = 0$. Thus
$(x_j)_{j=1}^{\infty}$ is a generic sequence of small points 
with respect to $\widehat{h}_{\bd{f}}$. 
\end{Example}

Now we prove the following theorem on equidistribution of small points 
on $\PP^1$ for $\bd{f}$. 
We remark that this theorem will also follow from 
Autissier \cite[Proposition~4.1.4 and Remarque 
in Introduction]{Au} (together with 
Proposition~\ref{prop:can:height:2}). We also remark that 
Baker--Hsia \cite{BH} proved equidistribution of small points on $\PP^1$ 
for polynomial maps $f$ over global fields satisfying the product formula. 

\begin{Theorem}
\label{thm:equidistribution}
Let $\bd{f} = (f_i)_{i=1}^{\infty}$ 
be a sequence of morphisms of $\PP^1$ over $K$ satisfying 
\textup{(C-i)} and \textup{(C-ii)}. 
Let $(x_j)_{j=1}^{\infty}$ a generic sequence of 
small points with respect to $\widehat{h}_{\bd{f}}$. Then 
$\frac{1}{\#G(x_j)}\sum_{y \in G(x_j)} \delta_y$ 
converges weakly to $T_{\bd{f}}$ as $j \to +\infty$.  
\end{Theorem}

\Proof
We adapt the proof of \cite{SUZ}. 
It suffices to show that, for any $C^{\infty}$ real-valued 
function $\varphi$ on $\PP^1(\CC)$, one has
\[
\lim_{j\to\infty} 
\frac{1}{\#G(x_j)}\sum_{y \in G(x_j)} \varphi(y) 
= T_{\bd{f}}(\varphi). 
\]
Let $\{(\calX_i, \overline{\calL}_i)\}_{i=0}^{\infty}$ 
be the adelic sequence associated 
with $\bd{f}$ in Proposition~\ref{prop:adelic:sequence}. 
We remark that $\calL_i$ is ample on $\calX_i$ because 
the normalization $f_{iU} \circ \cdots \circ f_{1U}: 
\PP^1_U \to \PP^1_U \hookrightarrow \PP^1_{O_K}$ is finite 
and $\calL_0$ is assumed to be ample on $\calX_0$.   
For $\lambda \in \RR$, we define the hermitian $\QQ$-line bundle 
$\overline{\calL}_i(\lambda\varphi)$ to be 
the $\QQ$-line bundle $\calL_i$ on $\calX_i$ equipped with 
the metrics $\{\Vert\cdot\Vert_{\sigma}^{\prime}\}_{\sigma}$ 
as follows: 
For all embeddings $\sigma$ except for $\sigma_0$ (resp. $\sigma_0$ and 
$\overline{\sigma_0}$) if $\sigma_0$ is real (resp. imaginary), 
the metric $\Vert\cdot\Vert_{\sigma}^{\prime}$ is the original 
metric $\Vert\cdot\Vert_{\sigma}$ of $\overline{\calL}_i$; 
For the embedding $\sigma_0$ (resp. the embeddings 
$\sigma_0$ and $\overline{\sigma_0}$), 
the metric $\Vert\cdot\Vert_{\sigma}^{\prime}$ 
is defined by 
\begin{align*}
& \Vert\cdot\Vert_{\sigma_0}^{\prime}  =
\Vert\cdot\Vert_{\sigma_0} \exp(- \lambda \varphi) 
\\
& \left(\text{resp. }
\Vert\cdot\Vert_{\sigma_0}^{\prime} = 
\Vert\cdot\Vert_{\sigma_0} \exp(- \lambda \varphi)
\;\;\;\text{and}\;\;\;
\Vert\cdot\Vert_{\overline{\sigma_0}}^{\prime} = 
\Vert\cdot\Vert_{\overline{\sigma_0}} \exp(- \lambda \varphi) \right).   
\end{align*}

\begin{Claim}
\label{claim:equidistribution}
$\displaystyle{
\liminf_{j \to\infty} h_{(\calX_i, \overline{\calL}_i(\lambda\varphi))} 
(x_j) \geq 
h_{(\calX_i, \overline{\calL}_i(\lambda\varphi))}(\PP^1)
}$. 
\end{Claim}

To see the claim, let $\varepsilon$ be any positive real number. 
Since $\deg(\calL_i\otimes_{O_K} K) = 1$, it follows from 
Theorem~\ref{thm:small:section} that there are a 
positive integer $n$ and a non-zero 
section $s \in H^0(\calX_i, \calL_i^{\otimes n})$ such that
\[
\Vert s\Vert_{\sigma, \sup}^{\prime} \leq \exp\left(
n \left(\varepsilon - 
\frac{\adeg\left(\acherncl_1(\overline{\calL}_i(\lambda\varphi))^2\right)}%
{2 [K:\QQ]}
\right) \right)
\]
for all embeddings $\sigma: K \hookrightarrow \CC$. 

Let $K_j$ be an extension field of $K$ such that 
$x_j$ is defined, and $O_{K_j}$ the ring of integers of $K_j$. 
Let $\Delta_{x_j}$ be the Zariski closure of $x_j$ in 
$\calX_i \times_{\Spec(O_K)} \Spec(O_{K_j})$. 
Let $p_j: \calX_i \times_{\Spec(O_K)} \Spec(O_{K_j}) \to \calX_i$ 
be the natural morphism. 
Since $(x_j)_{j=0}^{\infty}$ is a generic sequence, 
we have $\Delta_{x_j}  \not\subseteq 
\Supp(\zero(p_j^* s))$ for all $j \gg 1$. 
Then, for all $j \gg 1$, we get 
\begin{align*}
h_{(\calX_i, \overline{\calL}_i(\lambda\varphi))}(x_j) 
& = 
\frac{1}{[K_j:\QQ]} \adeg\left(\acherncl_1 
\left(\rest{p_j^* \overline{\calL}_i(\lambda\varphi)}{\Delta_{x_j}}
\right)\right) \\
& \geq  
\frac{1}{[K_j:\QQ]} 
\frac{1}{n} 
[K_j : K] 
\left(\sum_{\sigma: K \hookrightarrow \CC} 
- \log \Vert s\Vert_{\sigma, \sup}^{\prime}\right) \\
& \geq 
\frac{1}{[K_j:\QQ]}\frac{1}{n} 
[K_j : K] [K:\QQ] n 
\left(\frac{\adeg\left(\acherncl_1(\overline{\calL}_i(\lambda\varphi))^2\right)}%
{2 [K:\QQ]}
- \varepsilon \right)
= 
h_{(\calX_i, \overline{\calL_i}(\lambda\varphi))}(\PP^1) 
- \varepsilon. 
\end{align*}
Hence $\liminf_{j \to\infty} h_{(\calX_i, \overline{\calL}_i(\lambda\varphi))} 
(x_j) \geq 
h_{(\calX_i, \overline{\calL}_i(\lambda\varphi))}(\PP^1) - 
\varepsilon$. 
Since $\varepsilon$ is arbitrary, we obtain the claim. 

Now we compute $h_{(\calX_i, \overline{\calL}_i(\lambda\varphi))}(x_j)$ 
and $h_{(\calX_i, \overline{\calL}_i(\lambda\varphi))}(X)$. 
We set $r = 1$ if $\sigma_0$ is real and $r=2$ if 
$\sigma_0$ is imaginary.  
We have 
\begin{equation}
\label{eqn:equidistribution:1}
h_{(\calX_i, \overline{\calL}_i(\lambda\varphi))}(x_j)
= 
h_{(\calX_i, \overline{\calL}_i)}(x_j) 
+ 
\frac{\lambda r}{[K:\QQ]}  
\frac{1}{\#G(x_j)}\sum_{y \in G(x_j)} \varphi(y). 
\end{equation}

We also have 
\begin{align}
\label{eqn:equidistribution:2}
& h_{(\calX_i, \overline{\calL}_i(\lambda\varphi))}(\PP^1)
 = 
\frac{1}{2[K:\QQ]} \adeg\left(
\acherncl_1\left(\overline{\calL}_i(\lambda\varphi)\right)^2 \right) \\
& \qquad = \notag
\frac{1}{2[K:\QQ]} \adeg\left(
\acherncl_1\left(\overline{\calL}_i\right)^2 \right) 
+ 
\frac{\lambda r}{[K:\QQ]}  
\int_{\PP^1(\CC)} \varphi\,  
\cherncl_1(\overline{\calL}_i\otimes_{K^{\sigma_0}}\CC)
+
\frac{\lambda^2 r}{[K:\QQ]}  
\int_{\PP^1(\CC)} \varphi\, dd^c(\varphi) \\
&  \qquad = \notag
h_{(\calX_i, \overline{\calL}_i)}(\PP^1)
+ 
\frac{\lambda r}{[K:\QQ]}  
\int_{\PP^1(\CC)} \varphi\,  
\cherncl_1(\overline{\calL}_i\otimes_{K^{\sigma_0}}\CC)
+
\frac{\lambda^2 r}{[K:\QQ]}  
\int_{\PP^1(\CC)} \varphi\, dd^c(\varphi).
\end{align}
It follows from Claim~\ref{claim:equidistribution}, \eqref{eqn:equidistribution:1} and \eqref{eqn:equidistribution:2} that 
\begin{multline*}
\liminf_{j\to\infty}\left(h_{(\calX_i, \overline{\calL}_i)}(x_j) 
+ 
\frac{\lambda r}{[K:\QQ]}  
\frac{1}{\#G(x_j)}\sum_{y \in G(x_j)} \varphi(y)\right) \\
\geq 
h_{(\calX_i, \overline{\calL}_i)}(\PP^1)
+ 
\frac{\lambda r}{[K:\QQ]}  
\int_{\PP^1(\CC)} \varphi\,  
\cherncl_1(\overline{\calL}_i\otimes_{K^{\sigma_0}}\CC)
+
\frac{\lambda^2 r}{[K:\QQ]}  
\int_{\PP^1(\CC)} \varphi\, dd^c(\varphi).
\end{multline*}
Now we let $i$ to the infinity. 
Since the convergence $\widehat{h}_{\OO_{\PP^1}(1), \bd{f}}(x) 
= \lim_{i\to\infty} h_{(\calX_i, \overline{\calL}_i)}(x)$ is 
uniform with respect to $x \in \PP^1(\overline{K})$ by 
Proposition~\ref{prop:can:height:2}, we have 
\[
\lim_{i\to\infty}\liminf_{j\to\infty} 
h_{(\calX_i, \overline{\calL}_i)}(x_j) 
= \liminf_{j\to\infty} \lim_{i\to\infty} 
h_{(\calX_i, \overline{\calL}_i)}(x_j) 
= \liminf_{j\to\infty} 
\widehat{h}_{\OO_{\PP^1}(1), \bd{f}}(x_j) = 0. 
\]
By Theorem~\ref{thm:can:height}(2)(iv) and 
Proposition~\ref{prop:can:height:2}, 
we get $\lim_{i\to\infty} h_{(\calX_i, \overline{\calL}_i)}(\PP^1) 
= \widehat{h}_{\OO_{\PP^1}, \bd{f}}(\PP^1) = 0$. 
Moreover, by Lemma~\ref{lemma:estimate:Green:function}, 
Proposition~\ref{prop:Chern:current} and its proof, 
$\cherncl_1(\overline{\calL}_i\otimes_{K^{\sigma_0}}\CC)$ 
converges weakly to $T_{\bd{f}}$. Hence we get 
\[
\liminf_{j\to\infty}
\frac{\lambda r}{[K:\QQ]}  
\frac{1}{\#G(x_j)}\sum_{y \in G(x_j)} \varphi(y) 
\geq 
\frac{\lambda r}{[K:\QQ]} T_{\bd{f}}(\varphi)   
+
\frac{\lambda^2 r}{[K:\QQ]}  
\int_{\PP^1(\CC)} \varphi\, dd^c(\varphi).
\]
Since $\lambda$ is an arbitrary real number, we find 
$\lim_{j\to\infty} 
\frac{1}{\#G(x_j)}\sum_{y \in G(x_j)} \varphi(y) 
= T_{\bd{f}}(\varphi)$. This completes the proof. 
\QED

\end{document}